\documentclass{article}

\newcommand{\Om}{\Omega}
\newcommand{\la}{\langle}
\newcommand{\ra}{\rangle}

\newenvironment{pf}{\noindent{\sc Proof}.\enspace}{\rule{2mm}{2mm}\medskip}

\newtheorem{theorem}{Theorem}[section]
\newtheorem{proposition}{Proposition}[section]
\newtheorem{lemma}{Lemma}[section]
\newtheorem{corollary}{Corollary}[section]

\newtheorem{remark}{Remark}[section]
\newtheorem{remarks}{Remark}[section]
\newtheorem{definition}{Definition}[section]
\newcommand{\be}{\begin{equation}}
\newcommand{\ee}{\end{equation}}
\newcommand{\teta}{\theta}
\newcommand{\om}{\omega}

\newcommand{\al}{\alpha}

\newcommand{\ov}{\overline}
\newcommand{\wtilde}{\widetilde}
\newcommand{\R}{{\bf R}}

\newcommand{\Z}{{\bf Z}}

\newcommand{\N}{{\bf N}}
\newcommand{\T}{{\bf T}}
\newcommand{\sign}{{\rm sign}}
\renewcommand{\a }{\alpha }
\renewcommand{\b }{\beta }
\newcommand{\s }{\sigma }
\newcommand{\ii }{{\rm i} }
\renewcommand{\d }{\delta }

\newcommand{\e }{\varepsilon }
\newcommand{\g }{\gamma}

\renewcommand{\l }{\lambda }

\newcommand{\vphi}{\varphi }

\renewcommand{\t }{\tau }

\newenvironment{Remark}{\begin{remark} \rm}{\rule{2mm}{2mm}\end{remark}}

\newcommand{\dps}{\displaystyle}

\begin{document}

\title{{\bf Cantor families of periodic solutions for 
completely resonant nonlinear wave equations}}

\author{Massimiliano Berti\footnote{SISSA, 
Via Beirut 2-4, 34014, Trieste, Italy, {\tt berti@sissa.it}.}, \ 
Philippe Bolle\footnote{D\'epartement de math\'ematiques, Universit\'e
d'Avignon, 33, rue Louis Pasteur, 84000 Avignon, France, 
{\tt philippe.bolle@univ-avignon.fr}.}}

\date{}
\maketitle

{\bf Abstract:}
We prove existence of small amplitude,   
$2\pi \slash \om$-periodic in time solutions of completely resonant 
nonlinear wave equations with Dirichlet boundary 
conditions, for any frequency $ \om $ 
belonging to a Cantor-like set of positive measure 
and for a new set of nonlinearities. 
The proof relies on a suitable Lyapunov-Schmidt decomposition
and a variant of the Nash-Moser Implicit 
Function Theorem. 
In spite of the complete resonance of the equation
we show that we can still reduce the problem to 
a {\it finite} dimensional bifurcation equation. Moreover, 
a new simple approach for the inversion of the linearized
operators required by the Nash-Moser scheme is developed.
It allows to deal also with  nonlinearities which are not odd and 
with finite spatial regularity.
\\[2mm]
\footnote{Supported by 
M.I.U.R. Variational Methods and Nonlinear
Differential Equations.}
\\[1mm]
2000AMS subject classification: 35L05, 37K50, 58E05.

\section{Introduction and main result}

We consider the {\it completely resonant} nonlinear wave equation
\be\label{eq:main}
\cases{ u_{tt} - u_{xx} + f (x, u ) =  0 \cr
u ( t, 0 )= u( t,  \pi )  = 0}
\ee
where the nonlinearity $ f ( x, u ) = a_p (x) u^p + O ( u^{p+1} ) $ with  
$ p \geq 2$  is analytic with respect to $ u $ 
but has only finite regularity with respect to $ x$. 
More precisely, we assume
\begin{itemize}\label{regul}
\item[{\bf (H)}] There is $ \rho > 0 $ such that $ \forall (x,u) \in (0,\pi)
  \times (-\rho,\rho)$, $f( x , u ) =  \sum_{k=p}^{\infty} a_k(x)
  u^k $, $ p \geq 2$, where $a_k(x) \in H^1((0, \pi), \R)$ and 
$ \sum_{k=p}^{\infty} ||a_k||_{H^1} r^k < \infty$ for any $r\in
  (0,\rho)$. 
\end{itemize}
We look for small amplitude,  
$ 2\pi \slash \om$-periodic in time solutions of equation 
(\ref{eq:main}) for all
 frequency $ \om $ close to 
$1$ and  in some set of {\it positive measure}. 
\\[1mm]
\indent 
Equation (\ref{eq:main}) is an infinite dimensional Hamiltonian
system
possessing an elliptic equilibrium at $ u = 0 $. 
The frequencies of the linear oscillations at $0$ 
are $\om_j = j $, $ \forall j = 1, 2, \ldots $,   
and satisfy {\it infinitely many resonance} relations. 
Any solution $v = \sum_{j \geq 1} a_j \cos ( j t + \teta_j ) \sin (jx ) $ 
of the linearized equation at $u=0$, 
\be\label{eq:lin}
\cases{ u_{tt} - u_{xx} =  0 \cr
u ( t, 0 )= u( t,  \pi )  = 0}
\ee
is $2 \pi$-periodic in time.
For this reason equation (\ref{eq:main}) is called a completely 
resonant Hamiltonian PDE. 
\\[1mm]
\indent
Existence of periodic solutions close to a completely resonant
elliptic equilibrium for finite dimensional 
Hamiltonian systems 
has been proved by Weinstein \cite{We}, Moser \cite{Mo} and 
Fadell-Rabinowitz \cite{FR}. The proofs are based on the classical 
Lyapunov-Schmidt
decomposition 
which splits the problem into two equations: the {\it range equation}, 
solved through the standard Implicit Function Theorem, and 
the {\it bifurcation equation}, solved via variational arguments.  
\\[1mm]
\indent
For proving existence of small amplitude periodic solutions 
of completely resonant Hamiltonian PDEs like (\ref{eq:main}) 
two main difficulties must be overcome: 
\begin{itemize} 
\item[($i$)] a 
 ``{\it small denominators}'' problem which arises when solving the
 range equation; 
\item[($ii$)]  the  
presence of an {\it infinite dimensional} bifurcation equation:
which  solutions $ v $ of the linearized equation (\ref{eq:lin})
can be continued to solutions 
of the nonlinear equation (\ref{eq:main})? 
\end{itemize} 
The appearance of the small denominators problem $(i)$ is easily explained:
the  eigenvalues of the operator  $ \partial_{tt} -
\partial_{xx} $ 
in the spaces 
of functions  $u(t,x)$, $2 \pi / \om$-periodic in time 
and such that, say,  $u(t,.) \in  H^1_0 (0, \pi ) $ for all $t$, are 
$  - \om^2 l^2 + j^2 $, $l \in {\bf Z} $,
$ j \geq 1 $. 
Therefore, for almost every $ \om \in {\bf  R} $, 
the eigenvalues  accumulate to $ 0 $. As a consequence,  
for most  $\om$,  
the  inverse operator of $ \partial_{tt} -
\partial_{xx} $ is unbounded and the 
standard Implicit Function Theorem is not applicable.
\\[1mm]
\indent
The first existence results for small amplitude periodic solutions 
of (\ref{eq:main}) 
have been obtained 
in\footnote{Actually \cite{LS} deals with
the case of periodic boundary conditions in $x$,
i.e. $u(t,x+2\pi) = u(t,x)$.} 
 \cite{LS} for the nonlinearity  $ f ( x, u ) = u^3 $, 
and in \cite{BP1} for $f(x,u)=u^3+ O(u^5)$, 
 imposing on the frequency $ \om $ the 
 ``strongly non-resonance'' condition
$ |\om l - j | \geq \gamma / l $, $ \forall l \neq j $. 
For $\gamma >0$ small enough, the frequencies $\om $ satisfying such 
condition accumulate to $ \om =1 $ 
but form a set ${\cal W}_\g $ of zero measure.
For such $ \om $'s the spectrum of $ \partial_{tt} - \partial_{xx} $ 
does not accumulate to 
$ 0 $ and so 
the small divisor problem ($i$) is by-passed.
Next, problem $(ii)$ 
is solved by means of the Implicit Function Theorem, observing that
 the $0^{th}$-order bifurcation equation (which is an approximation of
 the exact bifurcation equation)
possesses, for $ f (x,u) = u^3 $,  non-degenerate periodic solutions,
see \cite{BP3}.

In \cite{BB}-\cite{BB1}, for  
the same set ${\cal W}_\gamma $ 
of strongly non-resonant 
frequencies, existence and multiplicity of 
periodic solutions has been proved for any  nonlinearity $ f(u) $. 
The novelty of \cite{BB}-\cite{BB1} was to solve the bifurcation equation 
via a variational principle
at fixed frequency which, jointly with 
min-max arguments, enables to find solutions of (\ref{eq:main})
as critical points of the Lagrangian action
functional. 

Unlike \cite{BP1}-\cite{BB}-\cite{BB1}, a new feature of the results
of this paper is that the set of frequencies $ \om $ 
for which we prove existence of $2\pi \slash \om $-periodic in time
solutions of (\ref{eq:main}) has positive measure, actually has
full density at $ \om = 1$.
\\[1mm]
\indent
The existence of periodic solutions 
for a set of frequencies of positive measure 
has been proved in \cite{B2} in the case of periodic boundary conditions
in $x$ and for the nonlinearity 
$ f(x,u) = u^3 + \sum_{4 \leq j \leq d} a_j (x)u^j $
where the $a_j (x) $ are trigonometric cosine polynomials in $x$.
The nonlinear equation $u_{tt} - u_{xx} + u^3 =0$ 
possesses a continuum of small amplitude,  
analytic and non-degenerate 
periodic solutions 
in the form of travelling waves
$ u(t,x) = $ $ \d p_0 ( \om t + x) $ where $ \om^2 = 1 + \d^2 $
and $p_0 $ 
is a non-trivial $2 \pi $-periodic solution 
of the ordinary differential equation $ p_0'' = - p_0^3 $. 
With these properties at hand
the small divisors problem ($i$) is solved 
via a Nash-Moser Implicit function Theorem 
adapting the estimates of Craig-Wayne \cite{CW}. 

Recently, the  existence of periodic solutions of (\ref{eq:main}) 
for frequencies 
$\om$ in a set of positive measure  has been proved in 
\cite{GMP} using the Lindstedt series method
for odd analytic nonlinearities $ f ( u ) = a u^3 + O( u^5 )$
with $ a \neq 0 $.
The reason for which $ f(u) $ must be odd is that the solutions are
obtained as  analytic sine-series in $ x $, 
see comments at the end of the section.

We also quote the recent papers 
\cite{IPT}-\cite{IPT1} on the
standing wave problem for a perfect fluid under gravity and with 
infinite depth which leads to a nonlinear and 
completely resonant second order equation.
\\[1mm]
\indent
In this paper
we prove the  existence 
of periodic solutions of the completely resonant wave equation
(\ref{eq:main}) with Dirichlet boundary conditions
for a  set of 
frequencies $\om$'s of positive measure and  
with full density at $\om = 1$ and  
for a new set of nonlinearities $ f(x,u) $ satisfying ({\bf H}) (including 
for example $f(x,u) = u^2$); 
we do {\it not} require that $f(x,u)$ can be extended 
on $ (-\pi,\pi) \times \R $ to an analytic 
function $g(x,u)$ satisfying the oddness assumption
$ g ( - x, - u ) = - g( x, u ) $,  
and we assume only finite regularity in the spatial variable $ x $,
see {\it Theorem \ref{thm:main}}.
\\[1mm]
\indent
Let us describe accurately our result.
Normalizing the period to $ 2 \pi $, we look for solutions 
$ u(t,x) $, $2 \pi$-periodic in time,  of the equation 
\be\label{eq:freq}
\cases{
\om^2 u_{tt} - u_{xx} + f (x, u ) = 0 \cr
u(t,0)= u(t, \pi) = 0
}
\ee
in  the real Hilbert space
\begin{eqnarray*}
X_{\s,s} := \Big\{  u (t,x) = \sum_{l \in {\bf Z}} 
\exp{(\ii lt)} \ u_l (x)   & \Big| &  u_l  \in H^1_0 ((0,\pi), {\bf C}), \ \ 
{\ov{u_l} (x)} = u_{-l} (x) \ \forall l \in {\bf Z}, \\ 
& & {\rm and} \  \ 
||u||_{\s,s}^2 := \sum_{l \in {\bf Z}} \exp{(2\s |l|)} 
(l^{2s} + 1 ) ||u_l||^2_{H^1} < +\infty \Big\}. 
\end{eqnarray*}

For $ \s > 0, s\geq 0 $, the space $ X_{\s,s}$                   
is the space of all $2\pi$-periodic in time functions with 
values in $ H^1_0 (( 0, \pi ), {\bf R}) $, which 
have a bounded analytic extension 
in the complex strip $|{\rm Im} \ t| < \s $ 
with trace function on $|{\rm Im } \ t| = \s $ belonging to  
$ H^s({\bf T}, H_0^1 ((0,\pi), {\bf C}))  $.

Note that if $u \in X_{\s,s}$ is a solution of (\ref{eq:freq}) in a
weak sense then the map $ x \mapsto u_{xx}(t,x) = \om^2
u_{tt}(t,x) - f(x,u(t,x)) $ belongs to $ H^1_0(0,\pi) $ for all $ t \in \T $;
hence $u(t,.) \in H^3(0,\pi) \subset C^2([0,\pi])$ and it is easy to
justify that $u$ is a classical 
solution.   

For $2 s > 1$, 
$ X_{\s,s} $ is a Banach algebra with respect to multiplication
of functions, namely\footnote{The proof is as in 
\cite{Po3} recalling that $H^1_0 (( 0, \pi ), {\bf C})$
is a Banach algebra with respect to multiplication of functions.} 
$$
u_1, u_2 \in X_{\s,s} \ \Longrightarrow \ 
u_1 u_2 \in X_{\s,s} \qquad {\rm and} \qquad || u_1 u_2 ||_{\s,s} \leq
C ||u_1 ||_{\s,s} ||u_2 ||_{\s,s}. 
$$
The space of the solutions 
of the linear equation $ v_{tt} - v_{xx} = 0 $ that belong to 
$ X_{\s,s} $ is\footnote{$ V $ can also be  written as
$$
V := \Big\{ v(t,x) = \eta ( t + x ) - \eta ( t - x ) \ 
\Big| \  \eta  \in C^\om_{\s,s+1} ({\bf T}, {\bf R}) \ {\rm with} \  
\int_{\bf T} \eta = 0 \Big\}
$$
where 
$C^\om_{\s,s+1} ({\bf T}, {\bf R}) $ denotes 
the space of all $2\pi$-periodic functions with have a bounded
analytic extension in the complex strip $|{\rm Im} \ t| < \s $ 
with trace function on $|{\rm Im } \ t| = \s $ belonging to 
$ H^{s+1} ( {\bf T}, {\bf C})$.}
$$
V := \Big\{  v (t,x) = \sum_{l \in {\bf Z}} 
\exp{(\ii lt)} u_l \sin (lx) \  \Big|  \ u_l \in {\bf C}, \ 
\ov{u_l} = - u_{-l}, \   
\sum_{l \in {\bf Z}} \exp{(2\s |l|)} 
(l^{2s} + 1 )l^2 |u_l|^2 < +\infty \Big\}. 
$$
Let $\e := \dps \frac{\om^2-1}{2}$.
Instead of looking for solutions of (\ref{eq:freq})
in a shrinking neighborhood of $0$ it is a
convenient devise to perform  the rescaling $ u \to \d u $ 
with $ \d := |\e|^{1/p-1} $
(in most cases, see however subsection \ref{subsec:u^2}), obtaining
\be\label{eq:freqre}
\cases{
\om^2 u_{tt} - u_{xx} + \e g (\d, x, u ) = 0 \cr
u(t,0)= u(t, \pi) = 0}
\ee
where 
$$
g(\d,x, u) :=  s^* \frac{f(x,\d u) }{\d^p} = s^* \Big( a_p(x) u^p + 
\d a_{p+1} (x) u^{p+1}
+ \ldots \ \Big)
$$
and $s^* :=\sign(\e )$, 
namely $ s^* = 1$ if $\om > 1$ and $s^*=-1$ if $ \om < 1$.
 
The main result of this paper is:

\begin{theorem}\label{thm:main} 
Consider the completely resonant nonlinear wave 
equation (\ref{eq:main}) where the nonlinearity 
$f(x,u) =  a_p(x) u^p + O ( u^{p+1} ) $, $ p \geq 2 $, 
satisfies assumption {\bf (H)}. \\[1mm]
1) There exists an open 
set ${\cal A}_p $ in $H^1((0,\pi), {\bf R})$ such that,
for all $a_p \in {\cal A}_p$, there is $ \s > 0 $ and a $C^\infty$-curve $
[0, \d_0) \ni \delta \to u ( \delta )(t,x) \in X_{\sigma,s}$
with the following properties:
\begin{itemize}
\item  $(i)$ 
There exists $s^* \in \{-1,1\}$ and  a Cantor set ${\cal C}_{a_p}
\subset [0,\delta_0)$ satisfying 
\be\label{meas}
\lim_{\eta \to 0^{+}} \frac{{\rm meas} ({\cal C}_{a_p} 
\cap (0, \eta))}{\eta} = 1
\ee
such that, for all $ \delta \in
{\cal C}_{a_p}$, $ u ( \delta ) $ is a $2\pi / \om$-periodic in time 
classical solution of (\ref{eq:main}) with
$ \om = \sqrt{2 s^* \d^{p-1}+1}$;
\item  $(ii)$ 
$|| \wtilde{u}(\delta) - \delta u_0 ||_{\s,s } = O( \delta^2 )$ 
for some $ u_0 \in V \backslash \{0 \}$,
where $\wtilde{u}(\delta) (t,x)=u(\delta) (t/\om ,x)$. 
\end{itemize}
All $ a_3(x) \in H^1(0,\pi)$ such that $\langle a_3 
\rangle := (1\slash \pi ) \int_0^\pi a_3 (x) dx \neq 0 $
is in ${\cal A}_3$. Hence $(i)$ and $(ii)$ hold true
for any nonlinearity like $ f(x,u) = a_3 ( x ) u^3 
+ \sum_{k\geq 4} a_k(x) u^k$, $\langle a_3 
\rangle \neq 0$, with $ s^*= {\rm sign}(\langle a_3 
\rangle) $.
\\[2mm]
2) In the case $ f(x,u) = a_2 u^2 + \sum_{k\geq 4} a_k(x) u^k$, $a_2 \neq 0 $,
 conclusions $(i)$ and $(ii)$ still hold true with 
$\om=\sqrt{-2\d^2+1}$.

\end{theorem}

\begin{remark}
($i$) Since equation (\ref{eq:main}) is autonomous,  any 
non-trivial time periodic solution $ u (t,x) $ of (\ref{eq:main}) 
generates a circle of solutions, i.e.
for any $ \theta \in \R $, 
$ u_\theta := u( t - \theta,x)$ is a solution too. 

($ii$) Multiplicity of periodic solutions $ u ( t, x )$ of 
(\ref{eq:main}) with increasing norm could be found, as in \cite{BB1}.
\end{remark}

\noindent
{\Large{\bf Sketch of the proof.}} We require first,
for simplicity of exposition (see note \ref{note4}), that 
\be\label{ap}
\exists v \in V \quad {\rm such \ that} \quad 
\int_{\Om} a_p (x) v^{p+1}(t,x) \ dt dx \neq 0, \qquad 
\Om := {\bf T} \times (0, \pi).
\ee 
To fix the ideas, we assume there is $ v \in V $ such that 
$ \int_\Om a_p(x) v^{p+1} > 0 $ and 
 we look for periodic solutions with 
frequency $ \om  > 1 $, so that $ s^* = 1 $ and 
$ \om = \sqrt{2 \delta^{p-1} + 1 }$
(if $ \int_\Om a_p(x) v^{p+1} < 0 $ for some $ v \in V $, then 
we can look for solutions of frequency
$ \om  < 1 $, so that  $ s^* = - 1 $ and $\om = \sqrt{-2 \delta^{p-1}+1} $).

Condition (\ref{ap})
is verified for $ p $ odd, iff $ a(x) $ is not antisymmetric
w.r.t. to $ x = \pi / 2 $, and for $ p $ even, iff
$ a(x) $ is not symmetric w.r.t. to $ x = \pi / 2 $, see 
Lemma \ref{vq} in the Appendix.
\\[1mm]
\indent
In order to find solutions of (\ref{eq:freqre}), 
we try to implement the usual Lyapunov-Schmidt reduction 
according to the decomposition 
$X_{\s,s} = V \oplus W  $ where 
\be\label{spW}
W :=
\Big\{ w = \sum_{l\in \Z} \exp (\ii lt)\  w_l(x) \in X_{\s,s} \ | \
 w_{-l}= \ov{w_l} \ \ {\rm and} \ \ 
\int_0^\pi w_l (x) \sin (lx) \ dx = 0, \ \forall l \in \Z \ \Big\}. 
\ee
Looking for solutions 
$ u = v + w $ with $ v \in V$, $ w \in W $ and
we are led to solve the bifurcation equation 
(called the ($Q$)-equation)
and the range equation (called the ($P$)-equation)
\be\label{eqs1}
\cases{
- \Delta v = \Pi_V g(\d, x,v + w) \  \ \ \qquad  \qquad (Q) \cr 
L_\om w = \e \Pi_W g(\d,x, v + w)   \ \qquad  \qquad (P)}
\ee
where 
$$ 
 \Delta v := v_{xx} + v_{tt}, \qquad 
\qquad   
L_\om := - \om^2 \partial_{tt} + \partial_{xx} 
$$ 
and 
$\Pi_V : X_{\s,s} \to V$,  $\Pi_W : X_{\s,s} \to W $
denote the projectors respectively on $V$ and $ W $.
\\[1mm]
\indent
Since $ V $ is infinite dimensional
a difficulty arises 
in the application of the method of \cite{CW} in presence
of small divisors:
 if $ v \in V \cap X_{\s_0,s} $ then the solution
$ w(\d , v) $ of the range equation, obtained with any Nash-Moser 
iteration scheme will have a lower regularity, e.g.  
$ w(\d, v) \in X_{\s_0 \slash 2,s} $. Therefore
in solving next the bifurcation equation for $ v \in V $, 
the best estimate we can obtain is $v \in V \cap X_{\s_0 \slash 2,s+2} $,
which makes  the scheme incoherent. 
Moreover we have to ensure that the $0^{th}$-order
bifurcation equation\footnote{The right hand side 
$ \Pi_V ( a_p(x) v^p ) $ is not identically equal to $0$ in $V$
iff condition (\ref{ap}) holds.
If not verified, as for $ f(x,u) = u^2 $,  
the $0^{th}$-order non-trivial bifurcation equation 
will involve  higher 
order nonlinear terms, see \cite{BB} and subsection
\ref{subsec:u^2}. \label{note4}}, i.e.
the ($Q$)-equation for $ \d = 0 $, 
\be\label{eq:unpe}
- \Delta v = 
\Pi_V \Big( a_p(x) v^p \Big) 
\ee 
has solutions $ v \in V $ which are analytic, 
a necessary property to initiate an analytic Nash-Moser scheme
(in \cite{CW}-\cite{CW1} this problem does not arise since,
dealing with non-resonant or partially resonant Hamiltonian PDEs like
$ u_{tt} - u_{xx} + a_1 (x ) u = f(x,u) $,
the bifurcation equation is finite dimensional).
\\[1mm]
\indent
We overcome  this difficulty thanks to a reduction    
to a {\it finite dimensional} bifurcation equation (on a subspace  
of $V$ of dimension $N$ independent of $ \om $). This reduction can
be implemented,   
in spite of  the complete resonance of  equation (\ref{eq:main}),
thanks to the compactness of the operator $(-\Delta)^{-1}$. 
 
We  introduce the  decomposition
$ V = V_1 \oplus V_2 $
where 
$$
\cases{
V_1 := \Big\{ v\in V  \ | \ v(t,x) = 
\sum_{l = 1}^N 
\Big( \exp{(\ii lt)} u_l + \exp{(- \ii lt)} \ov{u_l} \Big)
\ \sin (l x), \ u_l \in {\bf C} \Big\}  \cr 
V_2 := \Big\{ v \in V \ | \ v(t,x) = 
\sum_{l \geq N+1} \Big( \exp{(\ii lt)} u_l + \exp{(- \ii lt)}  \ov{u_l}\Big) 
\sin (l x), \ u_l \in {\bf C} \Big\} 
}
$$
Setting $ v := v_1 + v_2 $,  
with 
$ v_1 \in V_1, v_2 \in V_2 $, (\ref{eqs1}) 
is equivalent to 
\be\label{eqs}
\cases{
- \Delta v_1 = \Pi_{V_1} g(\d, x,v_1 + v_2 + w) \ \qquad  \qquad (Q_1) \cr 
- \Delta v_2 = \Pi_{V_2} g(\d,x,v_1 + v_2 + w)  \ \qquad  \qquad (Q_2) \cr
L_\om w = \e \Pi_W g(\d,x, v_1 + v_2 + w)   \ \qquad  \qquad (P)}
\ee
where 
$\Pi_{V_i} : X_{\s,s} \to V_i$ ($ i = 1, 2 $), 
denote the orthogonal projectors 
on $V_i$ ($ i = 1, 2 $). We specify that all the norms $|| \ ||_{\s
  ,s}$ are equivalent on $V_1$. In the sequel $B(\rho , V_1)$ will 
denote $B(\rho , V_1) := $ $ \{v_1 \in V_1 \ | \ ||v_1||_{0,s} < \rho \}$.  
\\[1mm]
\indent
Our strategy to find solutions of system (\ref{eqs}) - and hence to prove 
Theorem \ref{thm:main}- is the following. We solve first
the $(Q_2)$-equation obtaining 
$ v_2 = v_2 (\d, v_1, w) \in V_2 \cap X_{\s , s }$
by a standard 
Implicit Function Theorem provided we have chosen  $N$ large enough 
and $0 < \s < \ov{\s} $ small enough - depending on the nonlinearity $f$ 
but {\it independent of $\d$}, see section \ref{sec:Q2}.

Next we solve  
the $(P)$-equation, obtaining $w = w( \d, v_1) \in W 
\cap X_{\s \slash 2 ,s}$  
by means of a Nash-Moser Implicit Function Theorem
for $(\d , v_1) $ belonging to some Cantor-like set of parameters,
see section \ref{sec:P}. 

A major role  
is played by the inversion of the {\it linearized operators}
obtained at any stage of the Nash-Moser iteration.
As usual, the main difficulty in controlling such inverse operators
is due to the fact the diagonal elements may be arbitrarily small.

Our approach -presented in section \ref{sec:lin}- 
is different from the one  in \cite{CW}-\cite{Bo1}-\cite{B3} 
which is based on the Fr\"olich-Spencer technique \cite{FS}. 
It allows to deal  with  
nonlinearities $f(x,u) $ with finite regularity in the
spatial variable $x$ and 
which are not the restriction to $(0,\pi) \times \R$ of
a smooth odd function. 

We first develop $ u(t, \cdot ) \in H^1_0 ((0,\pi), {\bf R}) $ 
in time-Fourier expansion only
and we distinguish the ``diagonal part'' $ D = $
diag$\{ D_k \}_{k \in {\bf Z}}$ of the  operator
that we want to invert. Next, 
using Sturm-Liouville
theory (see Lemma \ref{interm}), we diagonalize each $D_k $ 
in a suitable basis of  
$ H^1_0 ((0, \pi), {\bf R}) $
(close, but different from $(\sin jx)_{j\geq 1} $).
Assuming a ``first order Melnikov non resonance condition''
(Definition \ref{Deltap})  
we prove that its eigenvalues are polynomially bounded
away from $ 0 $ and so we invert $ D $ with sufficiently good estimates
(Corollary \ref{cor:d-1}). 
The presence of the ``off-diagonal'' 
Toepliz operators requires  to analyze
the ``small divisors'': for our method
it is sufficient  to prove that
the product of two ``small divisors'' is larger than a constant
if the corresponding ``singular sites'' are close enough, see
Lemma \ref{sme1}.

If the nonlinearity $ f ( x, u ) $ can be extended to an
analytic (in both variables) odd function
then the Dirichlet problem on $ [ 0, \pi] $ is actually equivalent to the 
$ 2\pi $-periodic problem within the space of all odd functions
and a natural configuration space to use is 
$ Y :=$ $\{ u(x) =$ $\sum_{j \geq 1 } u_j \sin (jx)$ $|$ 
$\sum_{j} \exp{(2 a j) } j^{2 \rho } |u_j|^2 <$ $+ \infty  \}$.
On the other hand, for (still analytic) non odd nonlinearities
$f$, it is not possible 
in general to find a smooth periodic solution $u(t,.)$ that 
belongs to $Y$ for all $t$. For example,
if $f(x,u)=u^2$ then it is easy to see (deriving twice the equation w.r.t. 
$x$ and using that $ u_{xx}(t,0)=0 $) that
any smooth solution $u$ must satisfy $-u_{xxxx}(t,0)+$$2 u_x^2(t,0)=$ $0$. 
Moreover,
if $u$ is not trivial then $u_x(.,0)$ is not identically $0$. 
Hence $u_{xxxx}(.,0)$
does not vanish everywhere, which implies that $u(t,.)$ is not in $Y$ for all 
$t$.
For these reasons we shall consider 
as configuration space $H^1_0 ((0,\pi), {\bf R}) $ (and the solutions
that we shall find
satisfy $u (t, \ . ) \in H^1_0 (0,\pi) \cap H^3 (0,\pi)$). 

Finally (section \ref{sec:Q1})  
we solve the finite dimensional $(Q_1)$-equation for a new set of 
nonlinearities : for these nonlinearities, we can define a smooth
curve $(\d \mapsto v_1(\d)\in V_1)$ such that, if $\d$ belongs to 
some ``large'' set, $(\d , v_1(\d))$ gives rise to an exact solution
of equation (\ref{eq:main}).   
\\[3mm]
{\bf Acknowledgements:}
Part of this paper was written when the second
author was visiting S.I.S.S.A. in Trieste.
\\[3mm]
{\bf Notations:}
$B(R;X)$ denotes the open  ball of radius $ R $ in the space $ X $ 
centered at $ 0 $. 
$f(z) = O(g(z))$ means that there is a universal constant $C$ 
such that $|f(z)| \leq C| g(z)|.$ $\overline{z} $ is the complex
conjugated of $ z $.

\section{Solution of the ($Q_2$)-equation}\label{sec:Q2}

The $0^{th}$-order bifurcation equation (\ref{eq:unpe}) 
is the Euler-Lagrange equation of the functional 
$\Phi_0 : V \to {\bf R}$, defined by 
\be\label{Phi0}
\Phi_0 (v) = \frac{||v||_{H^1}^2}{2} - 
\int_\Om a_p(x) \frac{v^{p+1}}{p+1}
\ dx dt, \qquad  \Om = {\bf T} \times (0,\pi ), 
\ee
where $||v||^2_{H^1}=\int_\Om v_t^2+v_x^2 \ dx \ dt$.  
By the Mountain-pass Theorem \cite{AR} and condition
(\ref{ap}) (recall that we assume $\int_\Om a_p(x) v^{p+1} > 0$
for some $v \in V $), the critical set
\be \label{defK0}
K_{0,c}:= \Big\{ v \in V \ | \ \Phi_0'(v) =0, \ \Phi_0 (v) \leq c \Big\}
\ee   
is  non-trivial ({\it i.e} $\neq \{ 0\}$) for $c>0$ large enough
and is compact for the $ H^1 $-topology, see \cite{BB}.

In fact, by a direct bootstrap argument, $K_{0,c}$ is a 
compact subset of 
$ V \cap H^k (\Om )$, for any $ k \geq 0 $. Therefore 
$K_{0,c} \subset  V \cap C^\infty  (\Om )$
(even if $ a_p (x) \in H^1((0, \pi), {\bf R}) $ only,  because the
projection $ \Pi_V $ has a  regularizing effect in the variable $ x $)
and   we have  
({\it a priori estimate})
$$
\sup_{v \in K_{0,c}} ||v||_{0,s} 
< R < +\infty . 
$$   
By the analyticity assumption $({\bf H})$ on the nonlinearity $ f $
and the Banach algebra property of $X_{\s,s}$,  
the Nemitsky operator $ X_{\s,s} \ni u \to 
g(\d,x,u) \in X_{\s,s} $ is  in  $ C^{\infty} ( {\cal U}_\d,
X_{\s,s}) $, where ${\cal U}_\d = \{ u  \in  X_{\s,s}
\ | \  |\delta| \ ||u||_{\s,s} \leq \wtilde{\delta} \}$ and
$\wtilde{\d}$ depends on the radius of convergence $ \rho $ of the 
power series that defines $f(x,u) $.  
 
\begin{lemma} \label{vdue}
{\bf (Solution of the ($Q_2$)-equation)}
There exists $ \ov{\s} > 0, N \in {\bf N}_+,\d_0>0$ such that,  
$ \forall 0 \leq \s < \ov{\s} $, $ \forall ||v_1||_{0,s} \leq 2R $,
$ \forall  || w ||_{\s,s} \leq 1 $, $\forall \d \in [0, \d_0)$, 
there exists a unique $ v_2 = v_2(\d,  v_1, w) \in X_{\s,s} $
with $ ||v_2 (\d, v_1, w)||_{\s,s} \leq 1 $
which solves the $(Q_2)$-equation. 

Moreover, for any $v\in K_{0,c}$, $v_2(0, \Pi_{V_1}v,0)=\Pi_{V_2}v$.

Furthermore $ v_2(\d, v_1, w) \in X_{\s,s+2} $, 
$v_2 (\cdot, \cdot , \cdot ) \in 
C^\infty \Big( [0,\d_0) \times B(2R;V_1) \times B(1; W \cap X_{\s,s} ), V_2 
\cap X_{\s,s+2} \Big)$ and $Dv_2$, $D^2 v_2$ are bounded on
$[0,\d_0) \times B(2R;V_1) \times B(1; W \cap X_{\s,s} ) $.

Finally, if in addition $||w||_{\s,s'} < + \infty $  
for some $s' \geq s$, then (provided $\d_0$ has been chosen small
enough) $||v_2 (\d,  v_1, w)||_{\s,s'+2}
\leq $ $ K(s', ||w||_{\s,s'} )$.
\end{lemma}

\begin{pf}
Fixed points of the nonlinear operator 
$ {\cal N}(\d,v_1,w, \cdot ) : V_2 \to V_2 $ defined by 
$$
{\cal N}(\d, v_1, w, v_2 ) := 
(- \Delta )^{-1} \Pi_{V_2} g ( \d, x, v_1 + w + v_2 )  
$$
are solutions of equation $(Q_2)$. For $ w \in W \cap X_{\s,s} $, 
$ v_2 \in V_2 \cap X_{\s,s} $ we conclude that 
${\cal N}(\d, v_1, w, v_2) \in V_2 \cap X_{\s,s+2}$ since
$g(\d,x, v_1 + w + v_2 ) \in X_{\s,s} $ 
and because of the {\it regularizing property} of the operator
$ (- \Delta )^{-1} \Pi_{V_2} : X_{\s,s} \to V_2 \cap X_{\s,s+2}$.

Let $ B := \{ v_2 \in V_2 \cap X_{\s,s} \ | \ ||v_2||_{\s,s} \leq 1 \} $.
We claim that there exists $ N \in \N $,  
$ \ov{\s} >0 $ and $\d_0>0$, such that 
for any $ 0 \leq \s < \ov{\s} $, $||v_1||_{0,s} \leq 2R $,
$ || w ||_{\s,s} \leq 1 $, $ \d \in [0, \d_0 )$ the operator
$ v_2 \to {\cal N}( \d, v_1, w, v_2 ) $
is a contraction in $ B $, more precisely
\begin{itemize}
\item ($i$)
$||v_2||_{\s,s} \leq 1 \Rightarrow 
||{\cal N}(\d, v_1, w, v_2 )||_{\s,s} \leq 1$;
\item ($ii$) $ \forall v_2, \wtilde{v}_2 \in B   \Rightarrow 
||{\cal N}(\d, v_1, w, v_2 )- {\cal N}(\d, v_1, w, \wtilde{v}_2 )
||_{\s,s} \leq (1 / 2) || v_2 - \wtilde{v}_2 ||_{\s,s}.$
\end{itemize}
Let us prove $(i)$. $ \forall u \in X_{\s,s}$,   
$ || (- \Delta)^{-1} \Pi_{V_2} u ||_{\s,s} \leq 
(C / (N+1)^2) ||u||_{\s,s}$ and so 
$ \forall || w ||_{\s,s} \leq 1$,  
$|| v_1 ||_{0,s} \leq 2R $, $ \d \in [0, \d_0) $, 
\begin{eqnarray*}
|| {\cal N}(\d, v_1, w, v_2 )||_{\s,s} & \leq & 
\frac{C}{(N+1)^2} \Big| \Big| g(\d,x, v_1 + v_2 + w) \Big| \Big|_{\s,s} 
\leq  \frac{C'}{(N+1)^2} 
\Big(|| v_1||^p_{\s,s} + ||v_2||^p_{\s,s} + ||w||^p_{\s,s} \Big)\\
& \leq & \frac{C'}{(N+1)^2} 
\Big( \exp{(\s p N )} || v_1||^p_{0,s} + ||v_2||^p_{\s,s} + 1 \Big)
\leq  \frac{C'}{(N+1)^2} 
\Big( (4 R)^p  + ||v_2||^p_{\s,s} + 1 \Big)
\end{eqnarray*}
for $\exp{(\s N)} \leq 2 $, where we have used that 
$ || v_1||_{\s,s} \leq \exp{(\s N)} || v_1||_{0,s} \leq
4 R $.
For $ N $ large enough (depending on $R$) and $\ov{\s}= \ln 2 / N $, 
$(i)$ follows. 
Property $(ii)$ can be proved similarly and the existence of a
unique solution $v_2(\d,v_1,w) \in V_2 \cap X_{\s,s} $ follows. 

Since $K_{0,c}$ is compact in $V \cap H^s(\Om)$, we may assume that
$N$ has been chosen so large that for all $v \in K_{0,c}$, 
$||\Pi_{V_2} v||_{0,s} \leq 1/2$. Any $v \in K_{0,c}$ solves 
(\ref{eq:unpe}), hence $\Pi_{V_2}v$ solves the $(Q_2)$-equation
associated with $(0,\Pi_{V_1}v,0)$.  Hence $\Pi_{V_2}v \in B$ and 
$ \Pi_{V_2}v =  {\cal N}(0,\Pi_{V_1}v,0, \Pi_{V_2}v)$, which implies that 
$\Pi_{V_2}v=v_2 (0, \Pi_{V_1}v , 0)$. 

Because the map  $(\d, v_1, w,v_2 ) \mapsto {\cal N}(\d, v_1, w,v_2 )$ 
is $ C^{\infty} $, by the Implicit function Theorem
$ v_2  : \{ (\d,v_1,w)  \ | \ 
\d \in [0, \d_0 ) , || v_1 ||_{0,s} \leq 2R, ||w||_{\s,s} \leq 1 \} \to
V_2 \cap X_{\s,s} $ is a $C^{\infty}$ map.
Finally, since $ ( - \Delta )^{-1} \Pi_{V_2} $ is a continuous linear
operator from $X_{\s,s}$ to $ V_2 \cap X_{\s,s+2} $ and 
$$
v_2 (\d, v_1, w) =  ( - \Delta )^{-1} \Pi_{V_2} \Big( g ( \d, x, v_1 + 
w + v_2 (\d, v_1, w) \Big),
$$
by the regularity of the Nemitsky operator induced by $ g $, 
 $ v_2 (\cdot, \cdot , \cdot ) \in 
C^\infty ([0,\d_0) \times B(2R;V_1) \times B(1; W \cap X_{\s,s}), V_2 
\cap X_{\s,s+2})$.  The estimates for the derivatives can be obtained similarly.

Finally, assume that $||w||_{\s,s'} < +\infty$ for some $s'\geq s$. We
have
$$
v_2(\d,v_1,w)=(-\Delta)^{-1} \Pi_{V_2} g(\delta,x,v_1+w+v_2(\d,v_1,w)).
$$
Using the regularizing properties of the operator $(-\Delta)^{-1}
\Pi_{V_2}$ and the fact that $||v_1||_{\s,r} < +\infty$ for all $r
\geq s$, we can derive
by a simple bootstrap argument a bound on $||v_2(\d,v_1,w)||_{\s,s'}$.
It is useful here to observe that\footnote{This inequality can be obtained
considering the extension of the maps to the complex strip of width $ \s $, 
using that $ H^1_0(0,\pi)$ is a Banach algebra, and the inequality
$|| u^l ||_{H^r, {\rm per}} \leq $ $ K^l  ||u||_\infty^{l-1} 
||u||_{ H^r, {\rm per } } $ for $ 2 \pi $-periodic functions.
The last inequality follows from  
$|| uv ||_{H^r, {\rm per}} \leq $ $C_r (||u||_\infty 
||v ||_{H^r, {\rm per}} +$ $||v||_\infty||u ||_{r, {\rm per}})$
which is related to the Gagliardo-Nirenberg inequalities.}
$|| a_l(x) ( u^l ) ||_{\s,r} \leq C^l ||a_l||_{H^1} ||u||_{\s,s}^{l-1}
||u||_{\s,r}$ for $ r \geq s $, $l \geq p$, so that 
if $u\in X_{\s,r}$ and $\delta ||u||_{\s,s}$ is
small enough, then $g(\delta, x,u) \in X_{\s,r}$.   
\end{pf}

\begin{remark}
Lemma \ref{vdue} implies, in particular, that any solution  
$v \in K_{0,c} $ of the $0^{th}$-order bifurcation equation 
(\ref{eq:unpe}) is not only in $ V \cap C^\infty ( \Om )$ but actually belongs 
to $ V \cap X_{\s,s+2} $ and therefore is analytic
in $t$ and hence in $ x $. 
\end{remark}

We conclude this section with a Lemma which is a standard consequence
of the Lyapunov-Schmidt decomposition dealing with 
variational equations. 

Let us define the ``reduced'' functional
$ \Psi_0 : B(2R,V_1)= \{ v_1 \in V_1 \ | \ ||v_1||_{0,s} < 2R \} \to {\bf R} $ by 
\be\label{psi0}
\Psi_0 (v_1) := \Phi_0 (v_1 + v_2 (0,v_1, 0)).
\ee

\begin{lemma}\label{q20}
If $\ov{v}_1 $ is a critical point of 
$ \Psi_0 $,  
then $ \ov{v} := \ov{v}_1 + v_2 (0, \ov{v}_1, 0) $ is a critical point of
$\Phi_0 : V \to {\bf R}$, i.e.  $\ov{v} \in V $ is a solution
of the $0^{th}$-order bifurcation equation (\ref{eq:unpe}).
Conversely, if $\ov{v}$ is a critical point of $\Phi_0$ of critical
value $\leq c$ then there is a critical point $\ov{v}_1 \in B(2R,V_1)$
of $\Psi_0$ such that $\ov{v}= \ov{v}_1 + v_2 ( 0, \ov{v}_1 , 0)$.
\end{lemma}

\begin{pf}
Since $v_2 (0,v_1, 0)$ solves the $(Q_2)$-equation 
(for $\d = 0$, $w=0$),
$ d \Phi_0 ( v_1 + v_2 (0,v_1, 0)  )[k] = 0 $,
$ \forall k \in V_2 $.
Moreover, since  $\forall v_1 \in V_1 $ 
$v_2 (0,v_1, 0) \in V_2 $,
$ \partial_{v_1} v_2 (0,v_1, 0)[h] \in V_2$,
$ \forall h \in V_1 $. 
Therefore 
\begin{eqnarray}
d \Psi_0 (v_1) [h] & = & d \Phi_0( v_1 + v_2 (0,v_1, 0)  )
\Big[ h + 
\partial_{v_1} v_2 (0,v_1, 0)[h] \Big] =   
d \Phi_0( v_1 + v_2 (0,v_1, 0)  )[h ] \nonumber \\ 
& = &
\int_{\Omega} \Big[ - \Delta v_1 - \Pi_{V_1} \Big( a_p(x) (v_1 +
v_2 (0,v_1, 0)) \Big)^p\Big] h \ dx \ dt.\label{0Q1}
\end{eqnarray}
Therefore if $\ov{v}_1 $ is a critical point of 
$ \Psi_0 $, then $\ov{v} := \ov{v}_1 +v_2 (0, \ov{v}_1, 0)$
is a solution of  equation (\ref{eq:unpe}).

Conversely, assume that  $\ov{v}$ is a critical point of 
$ \Phi_0 $, with $\Phi_0 (v) \leq c$ and let 
$\ov{v}_1 = \Pi_{V_1}\ov{v}$. By Lemma \ref{vdue}, 
$\Pi_{V_2}\ov{v}=v_2(0, \ov{v}_1,0)$ and it is clear that $\ov{v}_1$
is a critical point of $\Psi_0$.
\end{pf}

\section{Solution of the ($P$)-equation}\label{sec:P}

By the previous section we are reduced to solve the $(P)$-equation 
with $ v_2 = v_2 (\d , v_1, w) $, namely
\be\label{eqP}
L_\om w = \e \Pi_W \Gamma (\d, v_1, w)
\ee
where 
$$ \Gamma (\d, v_1, w)(t,x) := 
g \Big( \d, x , v_1 (t,x) + w(t,x) + v_2(\d,v_1,w)(t,x) \Big).
$$ 
\indent
The solution $w= w(\d, v_1)$ of the $(P)$-equation (\ref{eqP}) is obtained 
by means of a Nash-Moser Implicit Function Theorem
for $( \d, v_1) $ belonging to a Cantor-like set of parameters.
\\[1mm]
\indent
We consider the orthogonal splitting 
$ W = W^{(n)} \oplus W^{(n)\bot}$ where 
\be\label{Wp}
W^{(n)}= \Big\{ w \in W \ \Big| \ w = \sum_{l=-L_n}^{ L_n} \exp{(\ii lt)} 
\ w_l(x) 
\Big\} 
\quad {\rm and} \quad 
W^{(n)\bot} = \Big\{ w \in W \ \Big| \ w = \sum_{|l|> L_n} 
\exp{(\ii lt)} \ w_l(x) \Big\},
\ee
and $ L_n $ are integer numbers
(we will choose $ L_n  = L_0 2^n $ with $L_0 \in {\bf N} $ large enough). 
We denote by 
$$
P_n : W \to W^{(n)} \qquad \qquad {\rm and} \qquad \qquad   
P^{\bot}_n: W \to W^{(n)\bot}
$$
the orthogonal projectors
onto $W^{(n)}$ and $W^{(n)\bot}$.
\\[1mm]
\indent
The convergence of the recursive scheme is based on 
 properties $(P1)$-$(P2)$-$(P3)$ below. 

\begin{itemize}
\item {\bf (P1)} ({\bf Regularity}) 
$\Gamma (\cdot ,\cdot , \cdot , ) \in 
C^\infty \Big( 
[0,\d_0) \times B(2R;V_1) \times B(1; W \cap X_{\s,s} ), X_{\s,s} \Big)
$. Moreover, $\Gamma$, $D\Gamma$ and $D^2 \Gamma$ are bounded on 
$[0,\d_0) \times B(2R,V_1) \times B(1;W \cap X_{\s,s} )$.

\end{itemize}

($P1$) is a consequence of the $C^\infty$-regularity of the 
Nemistky operator induced by $ g ( \d, x, u ) $ on $X_{\s,s}$ 
and of the properties of the map
$ v_2$ stated in Lemma \ref{vdue}.

\begin{itemize}
\item {\bf (P2)} ({\bf Smoothing estimate}) 
$\forall$ $w  \in  W^{(n)\bot} \cap X_{\s,s}$ 
and $\forall $ $0\leq \s' \leq \s$,
$||w||_{\s',s} \leq \exp^{(-L_n (\s-\s'))} ||w||_{\s,s}$. 
\end{itemize} 
The next property $(P3)$ is an {\it invertibility property} of the linearized
operator ${\cal L}_n(\d,v_1,w): W^{(n)} \to W^{(n)}$ defined by
\be\label{Lp}
{\cal L}_n(\d,v_1,w) [h] := L_\om h - \e P_n \Pi_W D_w \Gamma (\d,v_1,w) [h]
\ee
where $w$ is the approximate solution obtained 
at a given step of the Nash-Moser iteration.

The invertibility of ${\cal L}_n(\d,v_1,w)$ is obtained 
excising the set of  parameters $ ( \d, v_1 ) $ for which $0$ is
an eigenvalue of ${\cal L}_n(\d,v_1, w) $. 
Moreover, in order to have bounds for
the norm of the inverse operator ${\cal L}_n^{-1}(\d,v_1, w)$
which are sufficiently good for the recursive scheme,
we also excise the parameters $(\d, v_1)$ for which 
the eigenvalues of ${\cal L}_n (\d,v_1,  w) $ are too small.
We prefix some definitions.

\begin{definition}\label{def:MV} {\bf (Mean value)} 
For $ \Om := {\bf T} \times (0,\pi) $ we define  
$$ M ( \delta , v_1 , w) :=  \frac{1}{|\Om |}
\int_\Om \partial_u g \Big( \delta ,x,  v_1 (t,x) + w(t,x) +
v_2(\d,v_1,w)(t,x) \Big) \ dx dt.$$
\end{definition}

\begin{definition}\label{wss}
We define for $1 < \t < 2 $ 
$$
[w]_{\s,s} := \inf \Big\{ \sum_{i=0}^q \frac{||h_i||_{\s_i,s}}{
(\s_i-\s)^{\frac{2(\t-1)}{\b}} } 
\  ;
\  q \geq 1, \ 
\ov{\s} \geq \s_i > \sigma , \ \ h_i \in W^{(i)}, 
\ w = \sum_{i=0}^q h_i   \Big\}
$$
where $\b : = \frac{2 - \t}{2}$ 
and we set $[w]_{\s,s} := \infty$ if the above set is empty. 
\end{definition}

\begin{definition}\label{Deltap}
{\bf (First order Melnikov non-resonance condition)}
Let $ 0 < \gamma < 1 $ and $ 1 < \tau < 2 $. We define (recall that 
$\om = \sqrt{2 \d^{p-1} +1} $ and $ \e = \delta^{p-1}$)
\begin{eqnarray*}
\Delta_n^{\g,\t} (v_1,w) &:= & 
\Big\{ \delta \in [0,\d_0) \  \Big| \  
| \om k -j| \geq \frac{\gamma}{(k+j)^{\tau}}, \ 
\Big| \om k -j -\e \frac{M(\delta , v_1 , w)}{2j}\Big| 
\geq \frac{\gamma}{(k+j)^{\tau}} \\
& & \ \ \forall k \in {\bf N}, \ j \geq 1, \ k \neq j, 
\ \frac{1}{3 |\e |} < k, \ k \leq L_n, \  j \leq 2L_n \Big\}.
\end{eqnarray*}  
\end{definition}

We claim that:

\begin{itemize}
\item {\bf (P3)} {\bf (Invertibility of ${ {\cal L}_n}$)} 
There exist positive constants $ \mu $, $\d_0$ and $ C $ such that, 
if $[w]_{\s,s} \leq \mu$, $||v_1||_{0,s} \leq 2 R $ and 
$\d \in \Delta_n^{\g,\t} (v_1,w)\cap [0,\d_0)$ for some
$0 < \g < 1$, $ 1 < \t < 2 $, then 
${\cal L}_n (\d,v_1,w)$ is
invertible and the inverse operator  
${\cal L}_n^{-1}(\d,v_1,w): W^{(n)} 
\to W^{(n)}$ satisfies 
\be\label{p3}
\Big| \Big| {\cal L}_n^{-1} (\d,v_1,w)[h] \Big| \Big|_{\s,s} \leq 
\frac{C}{\gamma} (L_n)^{\tau - 1} ||h||_{\s,s}.
\ee
\end{itemize}

Property $(P3)$ is the real core of the convergence proof
and where the analysis of the small divisors enters into play. 
Property $(P3)$ is  proved in section \ref{sec:lin}.

\subsection{The Nash-Moser scheme}

We are going to define recursively 
a sequence 
$\{ w_n \}_{n\geq 0} $ with $ w_n = w_n ( \d,v_1) \in W^{(n)} $, 
$w_0 = 0 $, defined on smaller and smaller sets of ``non-resonant''
parameters $(\d,v_1) $, 
$A_n \subseteq $ $ A_{n-1} \subseteq \ldots \subseteq 
A_1 \subseteq $ $ A_0 :=$ 
$\{ (\d,v_1) \ | \ \d \in [0, \d_0), \ ||v_1||_{0,s} \leq 2R \}$. 
The sequence $(w_n(\d,v_1))$ will converge  
to a solution $ w(\d,v_1)$ of the $(P)$-equation
(\ref{eqP}) for $(\d,v_1) \in A_{\infty} := \cap_{n\geq 1} A_n $.   
The main goal of the construction is to show that, 
at the end of the recurrence, 
the set of parameters  
$ A_{\infty} := \cap_{n\geq 1} A_n $
for which we have the solution $ w (\d,v_1) $, 
remains sufficiently large.
\\[1mm]
\indent
We define inductively the sequence $\{ w_n \}_{n \geq 0}$. 
Define the ``loss of analyticity'' $ \gamma_n $ by
$$
\gamma_n := \frac{\gamma_0}{n^2+1}, \qquad \s_0 = \ov{\s}, 
\qquad \s_{n+1} = \s_n - \gamma_n,
\qquad  \forall \ n \geq 0 , 
$$
where we choose $ \gamma_0 > 0 $ small such that the ``total
loss of analyticity'' $\sum_{n\geq 0 } \gamma_n = $ $ \gamma_0 
\sum_{n \geq 0} 1 / (n^2+1) \leq$ $ \ov{\s}/ 2 $, i.e.
$ \lim_{n \to +\infty }\s_n \geq \ov{\s}/ 2 > 0$.   
We also assume 
$$
L_n := L_0 2^n, \qquad  \forall \ n \geq 0,
$$
for some large integer $L_0 $ specified in the next proposition. 

\begin{proposition} \label{Indu} {\bf (Induction)}
Let $ w_0 = 0 $ and
$ A_0 :=$ $\{ (\d,v_1) \ | \ \d \in [0, \d_0), \ ||v_1||_{0,s} \leq 2R \}$.
There exists $\e_0:= \e_0 (\g, \tau)$, $L_0:= L_0 (\g, \t)> 0 $ 
such that $\forall |\e | \g^{-1} <\e_0 $,
 there exists a sequence 
$\{ w_n \}_{n \geq 0}$, $ w_n = w_n ( \d , v_1 ) \in W^{(n)} $, 
of solutions of the equation
$$
L_\om w_n - \e P_n \Pi_W \Gamma (\d, v_1 , w_n ) = 0,
\leqno{(P_n)} 
$$
defined for 
$ ( \d, v_1 ) \in A_n \subseteq $ $ A_{n-1} \subseteq \ldots \subseteq 
A_1 \subseteq $ $ A_0 $, with 
$ w_n (\d, v_1) = \sum_{i=0}^{n} h_i (\d, v_1)$, $h_0 := w_0 = 0$,   
$ h_i = h_i (\d, v_1) \in W^{(i)}$  
satisfying 
$ || h_i ||_{\s_i,s} \leq |\e | \gamma^{-1} 
\ \exp{(-\chi^i)} $ for some  $ 1 < \chi < 2 $, 
$ \forall i = 0, \ldots , n $.
\end{proposition}

\begin{pf} Fix some $\chi \in (1,2)$.
We assume $ \e_0 \g^{-1} >0 $ small enough such that
\be\label{smallser}
\frac{\e_0}{\gamma} 
\sum_{i \geq 0} \exp{(- \chi^i)} \Big( \frac{1 + i^2}{\gamma_0}\Big)^{\frac{2(\t-1)}{\b}}
 < \mu. 
\ee
The proof proceeds by induction. We recall that $\beta := (2-\tau) / 2$
and that $\mu$ is defined in property $(P3)$.

Suppose we have already defined a solution 
$ w_n = w_n ( \d , v_1 ) \in W^{(n)} $ of  equation ($P_n$) 
satisfying the properties stated in the proposition.
We want to define 
$$
w_{n+1} = w_{n+1}(\d,v_1) := w_n (\d,v_1)+ h_{n+1} (\d,v_1), 
\qquad  h_{n+1} (\d,v_1) \in W^{(n+1)}
$$ 
as an {\it exact} solution of the equation
$$
L_\om w_{n+1} - \e P_{n+1} \Pi_W \Gamma (\d, v_1 , w_{n+1} ) = 0.
\leqno{(P_{n+1})}
$$
In order to find a solution $w_{n+1} = w_n + h_{n+1}$ 
of  equation ($P_{n+1}$) we write, for $ h \in W^{(n+1)} $,
\begin{eqnarray} 
\label{eqfin}
L_\om (w_n + h ) - \e P_{n+1} \Pi_W \Gamma(\d, v_1 , w_n + h ) 
& = & L_\om w_n - \e P_{n+1} \Pi_W \Gamma(\d, v_1 , w_n ) \\
& + & L_\om h - \e P_{n+1} \Pi_W D_w \Gamma(\d,  v_1 , w_n )[h] 
+ R ( h ) \nonumber \\
& = & r_n + {\cal L}_{n+1}( \d, v_1, w_n) [h ] + R ( h ),
\end{eqnarray}
where, since $ w_n $ solves  equation $(P_n)$, 
$$
\cases{r_n := L_\om w_n - \e P_{n+1} \Pi_W \Gamma(\d, v_1 , w_n ) = 
- \e P_n^\bot P_{n+1} \Pi_W \Gamma(\d, v_1 , w_n ) \in W^{(n+1)} \cr \\ 
R ( h ) := \
- \e P_{n+1} \Pi_W \Big( \Gamma(\d,  v_1 , w_n + h ) - 
\Gamma(\d,  v_1 , w_n ) - D_w \Gamma(\d,  v_1 , w_n )[h]\Big)}
$$
satisfy, by properties $(P1)$ and  $(P2)$, 
\be\label{smallrp}
||r_n||_{\s_{n+1},s} \leq |\e | \ C \exp{(- L_n \gamma_n)}
\Big| \Big|  P_{n+1} \Pi_W \Gamma(\d, v_1 , w_n ) \Big| \Big|_{\s_n,s} \leq
|\e | \ C'\exp{(- L_n \gamma_n)}
\ee
and, by property $(P1)$, $ \forall h, h' \in W^{(n+1)} $ 
with $|| h ||_{\s_{n+1},s}$, $|| h' ||_{\s_{n+1},s} $ small enough
\be\label{RP1}
\cases{ || R ( h )||_{\s_{n+1},s} \leq C|\e | \  || h ||_{\s_{n+1},s}^2 \cr
|| R ( h)- R(h') ||_{\s_{n+1},s} \leq 
C|\e | \ (|| h ||_{\s_{n+1},s} + || h' ||_{\s_{n+1},s}  ) \ || h - h' ||_{\s_{n+1},s}.}
\ee
Since
$|| h_i ||_{\s_i,s} \leq $ $ |\e | \g^{-1} \ \exp{(-\chi^i)}$, 
$ \forall i = 0, \ldots , n $,
by (\ref{smallser}),
\be
[w_n]_{\s_{n+1},s} \leq \sum_{i=0}^{n} 
\frac{|| h_i ||_{\s_i,s}}{(\s_i - \s_{n+1})^{ \frac{2(\t-1)}{\b}}} 
\leq \frac{|\e |}{\gamma} \ \sum_{i=0}^{n} 
\frac{ \exp{(- \chi^i) }}{\gamma_i^{2(\t-1) \slash \b}}\leq
\frac{|\e |}{\gamma}  \ \sum_{i \geq 0} 
 \exp{(- \chi^i) } \Big( \frac{1+ i^2}{\gamma_0} \Big)^{ \frac{2(\t-1)}{\b} }  < \mu.  
\ee
Hence by property $(P3)$, the linear operator ${\cal L}_{n+1}( \d, v_1, w_n): 
W^{(n+1)} \to W^{(n+1)}$ is invertible for 
$ ( \d, v_1 ) $ restricted to the set of parameters
\be
A_{n+1} := \Big\{ ( \d, v_1) \in A_n \ | \ 
\d \in \Delta_{n+1}^{\g,\t} (v_1, w_n) \Big\} \subseteq A_n.
\ee
Moreover 
\be\label{invern+1}
\Big| \Big| {\cal L}_{n+1}(\delta,v_1,w_n)^{-1} \Big| \Big|_{\s_{n+1},s}
\leq \frac{C}{\gamma} (L_{n+1})^{\tau-1}, 
\qquad \forall (\d, v_1) \in A_{n+1}. 
\ee 
By (\ref{eqfin}) equation $(P_n)$ is equivalent to the fixed point problem
\be  \label{fixedp}
w_{n+1}=w_n+h, \qquad {\cal G}(\d, v_1, w_n, h ) = h ,
\ee
for   the nonlinear operator 
${\cal G}(\d, v_1, w_n, \cdot ): W^{(n+1)} \to W^{(n+1)}$, defined
by
$$
{\cal G}(\d, v_1, w_n, h ):= - {\cal L}_{n+1}( \d, v_1, w_n)^{-1} \Big( r_n +  
R( h ) \Big).
$$

To complete the proof of the proposition we need 
the following Lemma.

\begin{lemma} {\bf (Contraction)}
$ { \cal G }( \d, v_1, w_n, \cdot )$ maps the ball
$ B(\rho_{n+1};W^{(n+1)}) $
$ := $ $ \{ w \in W^{(n+1)} 
\ | \ ||w||_{\s_{n+1},s} \leq 
\rho_{n+1} \} $ of radius
$\rho_{n+1} := |\e |  \g^{-1} \ \exp{(- \chi^{n+1})}$ into itself   
and is a contraction in this ball. 
\end{lemma}

\begin{pf}
By (\ref{invern+1}), (\ref{smallrp}) and (\ref{RP1}), 
\begin{eqnarray}\label{sm1}
|| {\cal G}( \d, v_1, w_n, h )||_{\s_{n+1},s} & = &
\Big| \Big| {\cal L}_{n+1}( \d, v_1, w_n)^{-1} 
\Big( r_n + R ( h ) \Big) \Big| \Big|_{\s_{n+1},s}\nonumber \\
\nonumber
&\leq & 
\frac{C}{\gamma}
(L_{n+1})^{\tau-1} \Big( || r_n ||_{\s_{n+1},s} + ||R ( h )||_{\s_{n+1},s}
\Big) \\ 
& \leq &
\frac{C'}{\gamma} (L_{n+1})^{\tau-1} \Big(  |\e | \ \exp{(-L_n \gamma_n)} + 
|\e| \ ||h ||_{\s_{n+1},s}^2 \Big).
\end{eqnarray}
By (\ref{sm1}), if $|| h ||_{\s_{n+1},s} \leq \rho_{n+1}$
then $|| {\cal G}( \d, v_1, w_n, h )||_{\s_{n+1},s} \leq $
$ C' (L_{n+1})^{\tau - 1}\g^{-1}|\e |  (  \exp{(-L_n \gamma_n)} +$  
$ \rho_{n+1}^2 ) \leq $ $ \rho_{n+1} $, provided that 
\be\label{contrp}
\frac{C'}{\gamma} 
(L_{n+1})^{\tau-1}   |\e | \ \exp{(-L_n \gamma_n)} \leq \frac{\rho_{n+1}}{2}
\qquad \ {\rm and} \qquad  \ 
\frac{C'}{\gamma} 
(L_{n+1})^{\tau-1} |\e | \rho_{n+1} \leq \frac{1}{2}.
\ee
It is easy to check that for $L_{n+1} := L_0 2^{n+1} $ and 
$\rho_{n+1} := |\e | \g^{-1} \ \exp{(- \chi^{n+1})}$ both inequalities in 
(\ref{contrp}) are satisfied $\forall  n \geq 0 $,   
choosing $ L_0 $ large enough and $ | \e | \g^{-1}  \leq \e_0  $ small enough. 
With similar estimates, using (\ref{RP1}), we can prove that 
$ \forall h, h' \in B(\rho_{n+1}; W^{(n+1)})$,    
$ || {\cal G}( \d, v_1, w_n, h') - {\cal G}( \d, v_1, w_n, h) 
||_{\s_{n+1},s}
\leq$ $(1/2) || h - h' ||_{\s_{n+1},s}$ again for  
$ L_0 $ large enough and $ | \e |\g^{-1} \leq \e_0  $ small enough, uniformly 
in $ n $, and we conclude that 
${\cal G}( \d, v_1, w_n, \cdot )$ is a contraction on $ B(\rho_{n+1};
W^{(n+1)} ) $. 
\end{pf}

By the standard Contraction Mapping Theorem we deduce the 
existence, for $ L_0 $ large enough and $ \e $ small enough, 
of a unique $ h_{n+1} \in W^{(n+1)}$ solving 
(\ref{fixedp}) and satisfying
$$
|| h_{n+1} ||_{\s_{n+1},s} \leq \rho_{n+1} = 
\frac{|\e |}{\g} \ \exp{(- \chi^{n+1})}. 
$$
Summarizing,  
$  w_{n+1} ( \d , v_1 )= w_{n} ( \d , v_1 ) + h_{n+1} ( \d , v_1 ) $ is 
a solution in $W^{(n+1)}$ of equation ($P_{n+1}$), defined for 
$ ( \d, v_1 ) \in A_{n+1} \subseteq $ $ A_{n} \subseteq \ldots \subseteq$  
$ A_1 \subseteq $ $ A_0 $, and 
$ w_{n+1} (\d, v_1) = \sum_{i=0}^{n+1} h_i (\d, v_1)$,  
$ h_i = h_i (\d, v_1) \in W^{(i)}$ 
satisfying 
$ || h_i ||_{\s_i,s} \leq |\e |\g^{-1} \ \exp{(-\chi^i)} $ 
for some  $ \chi \in (1, 2) $, 
$ \forall i =0, \ldots , n+1 $. 
\end{pf} 

\begin{remark}
A difference with respect to the usual ``quadratic'' Nash-Moser scheme, 
is that $ h_n ( \d, v_1) $ is found as an 
exact solution of  equation $(P_n)$. 
It appears to be more convenient  to prove the regularity of 
$ h_n (\d, v_1)$ 
with respect to the parameters $(\d, v_1)$, see Lemma \ref{lem:der}.
\end{remark}

\begin{corollary} {\bf (Solution of the ($P$)-equation)}
For $ ( \d, v_1) \in A_\infty := \cap_{n \geq 0} A_n $, 
 $ \sum_{i \geq 0} h_i (\d, v_1)  $  
converges  normally in $ X_{\ov{\s}/2, s} $ 
to a solution $ w (\d, v_1) \in W \cap X_{\ov{\s}/2,s} $ of equation (\ref{eqP})
with $|| w (\d, v_1) ||_{\ov{\s}/2,s} \leq C |\e | 
 \g^{-1} .$
\end{corollary}

\begin{pf}
By proposition \ref{Indu},
for $ ( \d, v_1) \in A_\infty := \cap_{n \geq 0} A_n $, 
$ \sum_{i=0}^\infty ||h_i( \d, v_1)||_{\ov{\s}/2 , s} < \infty$. Hence 
$ \sum_{i \geq 0} h_i (\d, v_1)  $  
converges  normally in $ X_{\ov{\s}/2, s} $ 
to some $ w (\d, v_1) \in W \cap X_{\ov{\s}/2,s} $, and we have 
\be
|| w (\d,v_1) ||_{\ov{\s}/2,s} \leq  
\sum_{i\geq 0} || h_i (\d, v_1)||_{\ov{\s}/2,s} \leq 
\sum_{i \geq 0} || h_i (\d, v_1)||_{\s_i,s} \leq
\sum_{i\geq 0} |\e |\g^{-1} \ \exp{(- \chi^i)} = O (|\e |  \g^{-1}).
\ee
Let us justify that      
$ L_\om w = \e \Pi_W \Gamma (\d,  v_1 , w)$.
Since $ w_n $ solves equation ($P_n$) ,  
\be \label{*}  
L_\om w_n =\e P_n \Pi_W \Gamma (\d,  v_1 , w_n) = 
\e  \Pi_W \Gamma (\d,  v_1 , w_n)- \e P_n^\bot \Pi_W \Gamma (\d,  v_1 , w_n).
\ee
We have 
$$
\Big\| P_n^\bot \Pi_W \Gamma (\d,  v_1 , w_n) \Big\|_{\ov{\s}/2,s} \leq 
 C  \exp(- L_n (\s_n - (\ov{\s} \slash 2))) \leq 
C  \exp(- \gamma_0 L_0 2^n / (n^2+1) ).
$$
Hence, by $(P1)$, the right hand side in (\ref{*}) converges in
$X_{\ov{\s}/2,s}$ to $\Gamma (\d,  v_1 , w)$.  Moreover $(L_\om w_n)
\to L_\om w$ in the sense of distributions. Hence 
$ L_\om w = \e \Pi_W \Gamma (\d,  v_1 , w)$. 
\end{pf}

Before proving the key property ($P3$) on the linearized operator
we prove a ``Whitney-differentiability'' property for   $ w ( \d, v_1 )$ 
extending $ w ( \cdot , \cdot ) $
in a smooth way on the whole $ A_0 $.   

For this, some bound on the derivatives of $h_n = w_n - w_{n-1} $ 
is required. 

\begin{lemma} \label{lem:der}
{\bf (Estimates for the derivatives of $ h_n $ and $ w_n $)}
For $|\e | \g^{-1} $ small enough and $ L_0 $ large enough,  
$ \forall n \geq 0 $, 
the function $(\d, v_1) \to  h_n (\delta , v_1) $ 
is in $ C^\infty (A_n, W^{(n)}) $ and
the $ k^{th}$-derivative $ D^k h_n (\delta , v_1) $
satisfies 
\be \label{pard}
\Big| \Big| D^k h_n (\delta , v_1) \Big| \Big|_{\s_n,s}
\leq \frac{|\e |}{\g} \  K_1(k ,\ov{\chi})  \exp(-\ov{\chi}^{n}),
\ee
for any $\ov{\chi} \in (0, \chi) $ 
and a suitable positive constant $K_1(k ,\ov{\chi}) $.

As a consequence, the function $(\d, v_1) \to 
w_n (\delta , v_1) = \sum_{i=1}^n h_n (\d, v_1) $ 
is in $ C^\infty (A_n, W^{(n)}) $ and
the $k^{th}$-derivative $D^k w_n (\delta , v_1) $
satisfies 
\be \label{estim3}
\Big| \Big| D^k w_n ( \delta , v_1 ) \Big| \Big|_{\s_n,s}
\leq  \frac{|\e |}{\g} \ K_2(k),
\ee
for  a suitable positive constant $ K_2 ( k )$. 
\end{lemma}

\begin{pf}
First, $ ||\partial_\l h_0||_{\s_0,s} = $ $ ||\partial_\l w_0||_{\s_0,s}= 0$
(we denote $ \l := (\d, v_1 ) $).
Next, assume, by induction, that $h_n = h_n (\d, v_1) $ is a 
$C^\infty$ map defined in $A_n$. We shall
prove that $ h_{n+1} = h_{n+1} (\d, v_1) $ is $C^\infty$ too. 
First recall that $ h_{n+1} = 
h_{n+1}(\d, v_1)  $ 
is defined, in Proposition \ref{Indu}, for $(\d, v_1) \in A_{n+1}$ 
as a solution in $ W^{(n+1)}$ of  equation $(P_{n+1})$, namely 
$$ 
U_{n+1} \Big( \delta , v_1 , h_{n+1}(\d, v_1) \Big) = 0, \leqno{(P_{n+1})}
$$ 
where
$$
U_{n+1} (\delta , v_1 ,h) := L_\om (w_n+h) -\e P_{n+1} \Pi_W 
\Gamma (\d,v_1,w_n+h). 
$$
The map $U_{n+1} : [0,\d_0) \times V_1 \times  W^{(n+1)} \to 
W^{(n+1)}$ is $C^{\infty}$, and
we claim that 
$ D_h U_{n+1} (\d , v_1 , h_{n+1}) = {\cal L}_{n+1}(\delta, v_1, w_{n+1})$
is invertible and that 
\be\label{invnear}
\Big| \Big| {\cal L}_{n+1}(\delta,v_1,w_{n+1})^{-1} 
\Big| \Big|_{\s_{n+1},s} \leq \frac{C'}{\gamma} (L_{n+1})^{\tau-1}. 
\ee
As a consequence, by the Implicit Function 
Theorem, the map $(\d, v_1) \mapsto  h_{n+1}(\d, v_1)$ is  in $C^\infty (A_{n+1}, W^{(n+1)})$.

Let us prove (\ref{invnear}). 
By Proposition \ref{Indu}, $ ||h_{n+1}||_{\s_{n+1},s} \leq C |\e | \g^{-1}
\exp (-\chi^{n+1})$. Hence, by $(P1)$, 
$$
\Big| \Big|{\cal L}_{n+1}(\delta,v_1,w_{n+1}) - 
{\cal L}_{n+1}(\delta,v_1,w_n) \Big| \Big|_{\s_{n+1},s}
\leq C | \e | \  ||h_{n+1}||_{\s_{n+1},s} \leq C 
\frac{\e^2}{\g} \exp (-\chi^{n+1}).
$$
Recalling (\ref{invern+1}), $ {\cal L}_{n+1}(\delta,v_1,w_n)$
is invertible and $||{\cal L}_{n+1}(\delta,v_1,w_n)^{-1}||_{\s_{n+1},s}
\leq (C \slash \gamma)  (L_{n+1})^{\tau-1}$.  
Hence, provided that $ \e \g^{-1} $ is small enough 
(note that $(L_{n+1})^{\tau-1}= 
(L_0 2^{n+1})^{\tau-1} << \exp(\chi^{n+1})$ for $n$
large),
${\cal L}_{n+1}(\delta,v_1,w_{n+1})$   
is invertible and (\ref{invnear}) holds.

We now prove in detail  estimate (\ref{pard}) for $ k = 1 $. 
Deriving  equation ($P_{n+1}$) with
respect to some coordinate $\lambda$ of $(\d,v_1) \in A_{n+1}$, 
we obtain 
$$
{\cal L}_{n+1}(\delta,v_1,w_{n+1}) \Big[ 
\partial_\lambda h_{n+1} (\d ,v_1) \Big] =
- (\partial_\lambda U_{n+1} ) \Big( \d,v_1,h_{n+1}(\d,v_1) \Big). \leqno(P'_{n+1})
$$
Hence, by (\ref{invnear}),  $ \partial_\lambda h_{n+1}$ satisfies the
estimate
\be \label{parh}
\Big| \Big| \partial_\lambda h_{n+1} \Big| \Big|_{\s_{n+1},s} 
\leq \frac{C}{\gamma} 
(L_{n+1})^{\tau-1} \Big| \Big| 
(\partial_\lambda U_{n+1} ) (\d,v_1,h_{n+1}) \Big| \Big|_{\s_{n+1},s}.
\ee
There holds
\be\label{partder}
(\partial_\lambda U_{n+1} ) (\d,v_1,h)=L_\om \partial_{\lambda} w_n -
\e P_{n+1} \Pi_W \Big[ (\partial_\lambda  \Gamma ) (\delta,v_1,w_n+h)+ 
\partial_w  \Gamma (\delta,v_1,w_n+h)[\partial_{\lambda} w_n]  \Big]
\ee
that we can write as
\be\label{piur}
(\partial_\lambda U_{n+1} ) (\d,v_1,h)=(\partial_\lambda U_{n+1} )
(\d,v_1,0) + r(\delta , v_1 ,h),
\ee
where $ r(\delta , v_1 ,h) := 
(\partial_\lambda U_{n+1} ) (\d,v_1,h) - (\partial_\lambda U_{n+1} )
(\d,v_1,0) $ satisfies, by $(P1)$,   
\be \label{estimr}
\Big| \Big| r(\delta , v_1 ,h) \Big| \Big|_{\s_{n+1},s} 
\leq C | \e | \  ||h||_{\s_{n+1},s}
\Big( 1+||\partial_\lambda w_n||_{\s_{n+1},s} \Big). 
\ee
Now, since $ w_n = w_n (\d, v_1) \in W^{(n)} $ solves   
equation $(P_n)$, 
$$ P_n U_{n+1}(\delta , v_1 ,0) = 
L_\om w_n -\e P_n \Pi_W \Gamma (\d,v_1,w_n)) = 0, \  \ 
 \forall (\d, v_1) \in A_n . $$ 
Hence differentiating w.r.t. $ \l $ we get  
$ P_n (\partial_{\lambda} U_{n+1})(\delta , v_1 ,0) = 0 $,
and so
\begin{eqnarray}
(\partial_{\lambda} U_{n+1})(\delta , v_1 ,0) & = & 
P_n^{\bot}(\partial_{\lambda} U_{n+1})(\delta , v_1, 0 ) \nonumber \\ 
& = & P_n^{\bot} L_\om \partial_\l w_n 
-\e P_n^{\bot} P_{n+1} \Pi_W \wtilde{\Gamma}(\d,v_1) \nonumber \\
&  = & -\e P_n^{\bot} P_{n+1} \Pi_W
\wtilde{\Gamma}(\d,v_1), \label{estin+1}
\end{eqnarray}
where, by (\ref{partder}), 
$ \wtilde{\Gamma}(\d,v_1) :=
(\partial_\lambda \Gamma) 
(\d ,v_1, w_n) + \partial_w \Gamma (\d,v_1,w_n)[\partial_\l w_n].$
By (\ref{estin+1}), $(P2)$, $(P1)$     
\begin{eqnarray} \label{estim2}
\Big| \Big| (\partial_{\lambda} U_{n+1})(\delta , v_1 ,0)
\Big| \Big|_{\s_{n+1},s}
&\leq & | \e | \  \exp (-L_n \g_n ) \Big| \Big| \Pi_W
\wtilde{\Gamma}(\d,v_1) \Big| \Big|_{\s_{n},s} \nonumber \\ 
& \leq & C | \e | \exp (- L_n \g_n )  
\Big( 1+ ||\partial_\l w_n||_{\s_n,s} \Big). 
\end{eqnarray}
Combining (\ref{parh}), (\ref{piur}), (\ref{estimr}), (\ref{estim2}) and 
the bound 
$ ||h_{n+1}||_{\s_{n+1},s} 
\leq | \e |\g^{-1} \exp{(- \chi^{n+1})} $, we get 
\begin{eqnarray*}
||\partial_\l h_{n+1}||_{\s_{n+1},s} & \leq & \frac{C}{\gamma} 
(L_{n+1})^{\tau-1} | \e |  
\Big( \frac{|\e |}{\g} \exp (-\chi^{n+1}) + \exp (- L_n \g_n )\Big) 
\Big(1+ ||\partial_\l w_n||_{\s_n,s}\Big) \nonumber \\ 
&\leq & C(\ov{\chi}) \frac{|\e  |}{\g} \ 
\exp (-\ov{\chi}^{n+1}) \Big( 1 + ||\partial_\l w_n||_{\s_n,s}\Big) \leq
C(\ov{\chi}) \frac{|\e  |}{\g}  \   
\exp (-\ov{\chi}^{n+1}) 
\Big( 1+ \sum_{i=0}^n || \partial_\l h_i ||_{\s_i,s} \Big)
\label{finn+1}
\end{eqnarray*}
for any $\ov{\chi} \in (1, \chi)$, $\e\g^{-1}$ small and $ L_0 $ large. 
Hence, setting $ a_n := || \partial_\l h_{n}||_{\s_{n},s} $  
we get 
$$
a_0  =0  \ \quad {\rm and} \quad a_{n+1} 
\leq  C(\ov{\chi}) \frac{|\e  |}{\g}  \  \exp
(-\ov{\chi}^{n+1}) \Big( 1+ a_0 + \ldots + a_n \Big)
$$
which implies 
$ || \partial_\l h_{n}||_{\s_{n},s} = a_{n} \leq  K(\ov{\chi}) | \e |\g^{-1} 
\ \exp(-\ov{\chi}^{n})$, $ \forall n \geq 0 $, 
for a suitable positive constant $ K ( \ov{\chi} )$. 
We can prove in the same way the general estimate (\ref{pard})
from which (\ref{estim3}) follows.
\end{pf}

Since, by (\ref{pard}), $ h_n ( \d, v_1 ) = O( \e \g^{-1} 
\exp{(- \ov{\chi}^n ))} $ with all its derivatives, and 
the ``non-resonant'' set $ A_n $ is obtained at each step  
deleting strips of size $ O(\g \slash L_n^\tau )$, 
we can define (by interpolation, say) 
a $ C^\infty $-extension $\wtilde w (\d, v_1)$ 
of $ w (\d, v_1) $ for all $(\d, v_1) \in A_0 $.
More precisely we can prove:

\begin{lemma} {\bf  (Smooth Extension 
$ \wtilde{w}$ of $ w $ on $A_0$)}\label{smoex}
Given $ \nu > 0 $, there exists a function 
$ \wtilde w  \in C^\infty ( A_0, W \cap X_{\ov{\s}\slash 2, s}) $
such that, if $ (\d, v_1) \in A_\infty := \cap_{n\geq 0} A_n $ and 
$ {\rm dist}( (\d, v_1), \partial A_n) \geq 2 \nu / L_n^3$, 
$ \forall n \geq 0 $, then $ \wtilde{w}( \d, v_1 )$ solves
the ($P$)-equation (\ref{eqP}). 

More precisely $ \wtilde{w}(\d, v_1) $ 
satisfies, $ \forall k \in {\bf N} $, 
\be \label{estim4}
\Big| \Big| D^k \wtilde{w} (\d, v_1) \Big| \Big|_{\ov{\s}/2 ,s} 
\leq \frac{|\e  |}{\g}  \frac{C(k)}{\nu^k},
\qquad \forall (\d, v_1)  \in A_0,
\ee 
for suitable constants $ C(k) > 0$. 
Moreover $ \wtilde{w} (\d, v_1) := \lim_{n \to +\infty} 
\wtilde{w}_n (\d, v_1) $
where $ \wtilde{w}_n  $  is in  
$ C^\infty$ $ (A_0,$ $ W^{(n)} ) $, and the sequence $(\wtilde{w}_n)$ 
 converges uniformly  in $ A_0 $ to $ \wtilde{w} $
for the norm $|| \ ||_{\ov{\s}/2, s}$, more precisely
\be \label{estim5}
\forall (\d,v_1) \in A_0, \quad \Big| \Big| \wtilde {w} (\d, v_1 ) - 
\wtilde{w}_n (\d, v_1 ) \Big| \Big|_{\ov{\s}/2 ,s} \leq 
C \frac{|\e  |}{\g} \ \exp(-\ov{\chi}^n).
\ee

\end{lemma}

\begin{pf}
First we endow $\R \times
V_1$ with the  Borelian positive measure defined by  $\mu(E)=m(L^{-1}(E))$, where 
$L$ is some automorphism from  $\R^{N+1}$ to $\R \times V_1$  and 
$m$ is the Lebesgue measure in $\R^{N+1}$. Let $\varphi : \R \times
V_1 \to  \R_+$ be a $C^{\infty} $-map supported in the open ball of radius $1$
centered at $0$ with $\int \varphi \ d\mu =1$. Let
$$
\wtilde{A}_n := \Big\{ \l = (\d, v_1 )  \in A_n \ | \ 
{\rm dist}(\l, \partial A_n) \geq \frac{\nu}{L_n^3} \Big\} \subset A_n,   
$$
where $ \nu $ is some small constant to be specified later. We define
$\varphi_n$, $\psi_n: \R \times V_1 \to  \R_+$ as
$$
\varphi_n (\l) := \Big( \frac{L_n^3}{\nu} \Big)^{N+1} \varphi \Big(
  \frac{L_n^3 }{\nu} \l \Big), \quad \quad 
\psi_n (\l) := \int_{\wtilde{A}_n} 
\varphi_n (\l - \eta) \ d\mu (\eta)
= \Big( \varphi_n * \delta_{\wtilde{A}_n} \Big) (\l),
$$
where $ \delta_{\wtilde{A}_n}$ is the characteristic
function of the set ${\wtilde{A}_n}$ namely, 
$\delta_{\wtilde{A}_n} (\l) := 1$ 
if $\l \in \wtilde{A}_n$ and
$\delta_{\wtilde{A}_n} (\l):=0$ if $\l \notin \wtilde{A}_n$. 
Clearly $ \varphi_n $ 
is a $C^\infty$ map supported in the open ball of radius 
$ \nu \slash L_n^3 $
centered at $0$ and $\int \varphi_n \ d\mu =1$. 
It follows that $ 0 \leq \psi_n (\l) \leq 1 $,
supp$\psi_n$$\subset$int$A_n$  and $\psi_n (\l)=1$ if 
$\l \in \wtilde{A}_n$ satisfies $  {\rm dist}(\l , \partial A_n) 
\geq 2\nu / L_n^3 $. 
Moreover $ \psi_n $ 
is $ C^\infty $  
and it is easy to check that  
$  |D^k \psi_n (\l)| \leq C(k)
(L_n^3 / \nu)^k$, $\forall k \in {\bf N} $, $ \forall \l \in \R \times V_1 $ 
 for a suitable positive constant $C(k)$. 

Now we can define $\wtilde{w}_n : A_0 \to W^{(n)}$ by 
$$
\wtilde{w}_0 (\l) := w_0(\l) = 0, \quad  
\quad  \wtilde{w}_{n+1} (\l) :=
\wtilde{w}_n (\l)+ \wtilde{h}_{n+1} (\l)
$$
where  $ \wtilde{h}_{n+1} (\l) := \psi_{n+1}(\l) h_{n+1}(\l) $ if
$\lambda \in A_{n+1}$ and $ \wtilde{h}_{n+1} (\l)=0$ if 
$\lambda \notin A_{n+1}$. 
(note that $h_{n+1}$ is $C^{\infty}$ because supp$\psi_{n+1} \subset $ int$A_{n+1}$). 
We define $\wtilde{w}_n \in C^{\infty} (A_0, W^{(n)} ) $ by  $ \wtilde{w}_n (\l) := $ $
\sum_{i=1}^n \wtilde{h}_i (\l )$.

By the bound $|D^k \psi_n (\l)| \leq C(k)
(L_n^3 / \nu)^k $ given above and (\ref{pard}) we obtain 
$$
\Big| \Big| D^k \wtilde{h}_{n+1}(\l) \Big| \Big|_{\s_{n+1},s} 
\leq \frac{|\e  |}{\g}  \ C(k,\ov{\chi})
\Big( \frac{L_{n+1}^3}{\nu} \Big)^k 
\exp(-\ov{\chi}^n), \qquad \forall k \in N,  \
\forall \l \in A_0, \
\forall n \geq 0.
$$
As a consequence, the sequence $(\wtilde{w}_n)$ (and all its derivatives) 
converges uniformally in $ A_0 $ 
for the  norm $|| \ ||_{\ov{\s}/2 ,s}$
on $ W $, to some function 
$ \wtilde{w}(\d, v_1) \in C^{\infty} 
(A_0, W \cap X_{\ov{\s}/2 ,s} ) $ which satisfies (\ref{estim4})
and (\ref{estim5}).

Note that if $ \l \notin A_\infty := \cap_{n\geq 0} A_n $ 
then the series $\wtilde{w}(\l) = \sum_{n \geq
  1} \wtilde{h}_n (\l)  $ is a finite sum. On the other hand, if 
$\l \in A_\infty$ and ${\rm dist}(\l, \partial A_n) 
\geq 2\nu / L_n^3$, $ \forall n \geq 0 $, then $ \wtilde{w}(\l)
= w( \l ) $ solves
the ($P$)-equation (\ref{eqP}). 
\end{pf}

We complete this part with the following Lemma which 
will be used in section \ref{sec:Q1}.
Define 
$$
B_n := \Big\{ (\delta , v_1) \in A_0  \  |  \  \delta \in
\Delta^{2\gamma, \t}_n (v_1,\wtilde{w}(\delta , v_1))   \Big\}
$$
where we have replaced $\gamma$ with $2\gamma$ in
the definition of $\Delta_n^{\g,\t}$, see Definition \ref{Deltap}. 

\begin{lemma}\label{Binfty} 
If $ \nu \g^{-1} > 0 $ and $ \e \g^{-1} $ are small enough, then  
$ B_n \subset \{ (\d, v_1) \in 
A_n \  | \  {\rm dist}( (\d, v_1) , \partial A_n) \geq 2\nu /
L_n^3 \}$,  $\forall n \geq 1 $ and hence,
if $ (\delta,v_1) \in B_\infty 
:= \cap_{n\geq 1} B_n $, then
$\wtilde{w} (\d, v_1)$ solves the $(P)$-equation (\ref{eqP}). 
\end{lemma}

\begin{pf}
This is a consequence of (\ref{estim4}), (\ref{estim5}) and the
previous Lemma
(we use that $L_n^{\tau}=o(\exp(\wtilde{\chi}^n))$ and
$L_n^\tau = o(L_n^2)$ as $n \to \infty$. 
\end{pf}

\begin{Remark}
Up to now, we have not justified that $A_\infty$, $B_\infty$ are not empty.  
It can be proved as in Proposition \ref{measure} that for any 
$v_1 \in B(2R,V_1)$ the set ${\cal B}_{v_1}:=$ $\{ \d \in [0,\d_0) \ | \
(\d , v_1) \in B_\infty \}$ satisfies $\lim_{\eta \to 0+} {\rm meas}((0,\eta)
\cap {\cal B}_{v_1})/ \eta =0$. Hence $A_\infty \neq \emptyset$. 
\end{Remark}

\section{Analysis of the 
linearized problem: proof of (P3)}\label{sec:lin}

We prove in this section the key property ($P3$) on the inversion of the linear
operator $ {\cal L}_n (\d,v_1,w): W^{(n)} \to W^{(n)}$ defined in 
(\ref{Lp}). Let
$$
a(t,x) :=  \partial_u g \Big( \d ,x, v_1(t,x)+w(t,x)+v_2(\d, v_1,w)(t,x)\Big)
$$
and define the linear operators $D $, $ M_1 $, $M_2 : W^{(n)} \to
W^{(n)}$ by
\be\label{defin}
\cases{ 
Dh := L_\om h -\e P_n \Pi_W ( a_0(x) \  h) \cr
M_1 h := \e P_n \Pi_W (\ov{a}(t,x)\  h)  \cr
M_2 h := \e P_n \Pi_W (a(t,x) \  \partial_w v_2 [h]) }
\ee
where 
$$
\cases{
a_0(x) := (1 \slash 2\pi )  \int_0^{2\pi} a(t,x) \ dt \cr 
\ov{a}(t,x) := a(t,x)-a_0(x).}
$$
\indent 
$ {\cal L}_n (\d,v,w) $ 
can be written as
\begin{eqnarray*}
{\cal L}_n ( \d , v, w ) [ h ] & := &  
L_\om h - \e P_n \Pi_W D_w \Gamma ( \d , v_1 , w ) [ h ] \\
& = & L_\om h - \e P_n \Pi_W \Big( 
\partial_u g (\d ,x, v_1+w +v_2(\d, v_1,w) ) \Big( h 
+ \partial_w v_2 (\d, v_1, w) [h] \Big) \Big) \\
& = & L_\om h - \e P_n \Pi_W \Big( a(t,x) \  h \Big) - 
\e P_n \Pi_W \Big( a(t,x) \ \partial_w v_2 (\d, v_1, w) [h] \Big) \\
&= & Dh - M_1h - M_2 h.
\end{eqnarray*}
First (Step $1$) we prove that,
assuming the ``first order Melnikov non-resonance condition
$ \d \in \Delta_n^{\g,\t} (v_1, w)$ (see Definition \ref{Deltap})
 the linear operator $ D $ is invertible, see Corollary \ref{cor:d-1}.
Next (Step $2$) we prove that 
$ M_1 $, $ M_2 $ are small enough with respect to 
$ D $, yielding the invertibility of the whole $ { \cal L}_n $.  
\\[2mm]
\indent
Through this section we shall use the notations 
$
F_k := \{ f \in  H^1_0  ((0,\pi); {\bf C}) \ | \ \int_0^\pi f(x) \sin (kx) \ dx =0 \}
$
whence the space $ W $, defined in (\ref{spW}), 
and its corresponding projector $\Pi_W : X_{\s,s} \to W $, 
are written, for any $ h = \sum_{k \in {\bf Z}} \exp{(\ii kt) } h_k $ 
$$
W = 
\Big\{ h \in X_{\s,s} \ |  \ h_k \in F_k \ \forall k \in {\bf Z}  \Big\}, 
\qquad \Pi_W  h (t,x) = \sum_{k \in {\bf Z}} \exp{(\ii kt)} (\pi_k h_k) (x),
$$ 
where $\pi_k :  H^1_0  ((0,\pi); {\bf C}) \to F_k$ is 
the $L^2$-orthogonal projector onto $ F_k
$ (note that $\pi_{-k}=\pi_k$ and $\pi_k \ov{u}= \ov{\pi_k u}$, hence 
$\pi_{-k} h_{-k}=\ov{\pi_k h_k}$). 
\\[3mm] 
\noindent
{\large {\sc {\bf Step 1: Inversion of ${D : W^{(n)} \to W^{(n)}}$}}}.  
\\[3mm]
\indent
In term of  time-Fourier series, $D$ is defined by 
$ \forall \ h \in W^{(n)} $, 
$$
(Dh)_k = D_k h_k \qquad \forall \  |k | \leq L_n,  
$$
where $ D_k : {\cal D}(D_k) \subset F_k \to F_k $ is the linear operator 
$$
D_k u = \om^2 k^2 u - S_k u \qquad {\rm and} \qquad
S_k u := - \partial_{xx} u + \e \pi_k (a_0(x) \  u).
$$
Note that  $ S_k = S_{-k} $. 

We now analyze the spectral properties of the Sturm-Liouville type
operator $ S_k $. 
We shall assume that
$ |\e | |a_0|_\infty < 1 $,  so that 
$$
\la u , v \ra_\e  := \int_0^\pi u_x \ov{v}_x +  \e a_0(x) u \ov{v} \ dx
$$
defines a scalar product on $  H^1_0  ((0,\pi); {\bf C}) $, hence on
$F_k$,  and 
its associated norm is
equivalent
to the standard $H^1$-norm. More precisely \footnote{$\forall u \in H^1_0 ( 0,\pi ) $,  
$\int_0^\pi u^2 (x) \ dx \leq \int_0^\pi u_x^2 (x) \ dx $
since the least eigenvalue of $- \partial_{xx} $
with Dirichelet B.C. on $(0,\pi)$ is
$ 1 $.} ,
\be\label{eq:ue1}
|| u ||_{\e} = ||u||_{H^1}  \Big( 1 + O( \e |a_0|_{\infty})\Big) \qquad 
\forall u \in F_k.
\ee

\begin{lemma} \label{interm}
{\bf (Sturm-Liouville)} 
The operator $ S_k : {\cal D}(S_k) \subset F_k \to F_k $ possesses
a $ \la \ , \ \ra_\e$-orthonormal basis $ ( v_{k,j})_{j\geq 1, j\neq |k|}$
of real eigenvectors with real 
eigenvalues $( \lambda_{k,j})_{j \geq 1, j \neq |k|}$,  
i.e. $ S_k v_{k,j}= \lambda_{k,j} v_{k,j }$, 
$\lambda_{k,j} \in {\bf R}$, and  $\l_{k,j} = \l_{-k,j}$,
$v_{-k,j}=v_{k,j}$.

Moreover, $ (v_{k,j})_{j\geq 1, j\neq |k|} $ is 
an orthogonal basis also for the $L^2$-scalar 
product in $F_k$. 
Defining $\varphi_{k,j}= v_{k,j}/ ||v_{k,j}||_{L^2}$,
 $\lambda_{k,j}$ and  $\vphi_{k,j}$ 
have the asymptotic expansion as $j \to +\infty $ 
\be\label{expa} 
\lambda_{k,j} = \lambda_{k,j}( \d, v_1, w) =
j^2 + \e M(\d, v_1,w)  + O \Big( 
\frac{\e ||a_0||_{H^1}}{j} \Big), \qquad
\Big| \vphi_{k,j} - \sqrt{ \frac{2}{\pi}} \sin (jx) \Big|_{L^2} = O \Big(
\frac{ \e |a_0|_\infty }{j} \Big)
\ee
where $M (\d, v_1,w)$, introduced in Definition \ref{def:MV},  
is the mean value of $ a_0 ( x ) $ on $ ( 0, \pi )$.
\end{lemma}

\begin{pf}
In the Appendix.
\end{pf}

By Lemma \ref{interm}, the linear operator 
$ D_k : {\cal D}(D_k) \subset F_k \to F_k $ possesses
a $\la \ , \ \ra_\e$-orthonormal basis $ ( v_{k,j})_{j\geq 1, j\neq |k|}$
of real 
eigenvectors with real eigenvalues $( \om^2 k^2 -\lambda_{k,j})_{j \geq 1, 
j \neq |k| }$.  
As a consequence we derive

\begin{corollary} \label{cor1} {\bf (Diagonalization of $D$)}
The operator $ D : W^{(n)} \to W^{(n)} $ is the diagonal operator 
{\rm diag}$\{ \om^2 k^2 - \lambda_{k,j} \}$ 
in the basis $ \{ \vphi_{0,j} \ ; \ j \geq 1 \} \cup ( \cup_{j,k\geq 1,
  j \neq k }  \{ \cos(kt) \varphi_{k,j} , \sin(kt) \varphi_{k,j} \}$
of $W^{(n)}$.
\end{corollary}

Noting that $\min_{|k | \leq L_n} |\om^2 k^2 -\lambda_{k,j}| \to
\infty$ as $j\to + \infty$, we deduce from Corollary \ref{cor1} that  
the  linear operator $ D : W^{(n)} \to W^{(n)} $ is invertible
iff all its eigenvalues $ \{ \om^2 k^2 - \lambda_{k,j} ( \d, v_1, w) $ 
$ \}_{|k| \leq L_n, j \geq 1, j \neq |k|} $ are different from zero. 
If this holds, we can define $ D^{-1} $ as well as 
$|D|^{-1/2} : W^{(n)} \to W^{(n)}$ by 
$$
W^{(n)} \ni h = \sum_{ k=- L_n}^{L_n} \exp{(\ii kt)} h_k \longrightarrow  
|D|^{-1/2} h := 
\sum_{ k= -L_n}^{L_n} \exp{(\ii kt)} |D_k|^{-1/2} h_k \in W^{(n)}
$$
where $|D_k|^{-1/2} : F_k \to F_k $ is the compact operator defined by
$$
|D_k|^{-1/2} \vphi_{k,j} := 
\frac{\vphi_{k,j}}{\sqrt{|\om^2 k^2 -\lambda_{k,j}|}}, 
\qquad \forall j \geq 1, \ j \neq | k |.
$$

The ``small divisor problem'' ($i$) is that
some of the eigenvalues of $D$, $\om^2 k^2 - \lambda_{k,j}$, 
can become arbitrarily small for $(k,j) \in {\bf Z}^2 $ sufficiently large 
and therefore
the norm of $  | D |^{-1/2} $ can become arbitrarily large as $n\to \infty$.

In order to quantify this phenomenon we define for all $k$
\be
\alpha_k := 
\min_{j \neq |k|} | \om^2 k^2 - \lambda_{k,j}|.
\ee
Note that $\a_{-k} = \a_k $. 
If $\forall |k | \leq L_n, \  \a_k \neq 0 $,  then
$ D $ is invertible and, since $(  \varphi_{k,j} )_{k \neq j}  $  is
an orthogonal basis for the $\la \ , \ \ra_\e $ scalar product,
$|| \ |D_k|^{-1/2} u \ ||_{\e} \leq $ $\alpha_k^{-1\slash 2} 
||u||_{\e} $. Hence, by (\ref{eq:ue1}), 
\be \label{esdk}
\Big| \Big| |D_k|^{-1/2} u \Big| \Big|_{H^1} \leq 
\frac{C}{\sqrt{\alpha_k}} ||u||_{H^1}, \qquad \forall k \in {\bf Z}. 
\ee
The condition  $\forall |k | \leq L_n, \  \a_k \neq 0   $, 
depends very sensitively on the parameters $(\d, v_1)$.
Assuming the ``first order Melnikov non-resonance
condition'' $ \d \in \Delta_n^{\g,\t} (v_1,w) $
(see Definition \ref{Deltap}), we obtain, in Lemma \ref{lowb}, a 
 lower bound of the form $c/ |k|^{\tau-1}$  for the moduli of the 
 eigenvalues of $ D_k $ and,
therefore, in Corollary
\ref{cor:d-1}, 
sufficiently good estimates for the inverse of $ D $.

\begin{lemma} {\bf (Lower bound for the eigenvalues of $D$)}\label{lowb}
If $ \delta \in \Delta_n^{\gamma,\t} (v_1,w) \cap [0, \d_0) $ and 
$\d_0  $ 
is small enough (depending on $\gamma$), then (recall that $ 1 < \tau < 2 $ )
\be\label{akl}
\alpha_k := \min_{j \geq 1, j \neq |k|} | \om^2 k^2 - \lambda_{k,j}| 
\geq \frac{\gamma}{|k|^{\tau-1}} > 0, \qquad \forall \ |k | \leq L_n.
\ee
\end{lemma}

\begin{pf}
Since $\a_{-k}= \a_k $ it is sufficient to consider $ k \geq 0 $.
By the asymptotic expansion (\ref{expa}) for the eigenvalues
$ \lambda_{k,j} $, using that $ || a_0 ||_{H^1}$, $ | M(\d, v_1,w) | \leq C $,
\begin{eqnarray}\label{estilb}
| \om^2 k^2 - \lambda_{k,j} | & = & \Big| \om^2 k^2 - j^2 - 
\e M(\d, v_1,w)  + O \Big( \frac{\e ||a_0||_{H^1}}{j} \Big) \Big| \nonumber \\
& = & \Big| \Big( \om k - \sqrt{j^2 + \e M(\d, v_1,w) } \Big) 
\Big( \om k + \sqrt{j^2 + \e M(\d, v_1,w) } \Big) + O \Big( \frac{|\e |}{j}
\Big) \Big|\nonumber \\
& \geq & \Big| \om k -j - \e \frac{M(\delta , v_1 , w)}{2j} + 
O \Big( \frac{\e^2 }{j^3} \Big) \Big| 
\om k - C \frac{|\e | }{j} \nonumber \\
& \geq & \Big| \om k -j - \e \frac{M(\delta , v_1 , w)}{2j} 
\Big| \ \frac{k}{2} 
- C' \Big( \frac{\e^2 k }{j^3} + \frac{|\e |}{j} \Big)
\geq \frac{c \gamma}{k^{\tau-1}} -  C \Big( \frac{\e^2 k }{j^3} 
+ \frac{|\e |}{j}\Big),
\end{eqnarray}
since $\d \in \Delta_n^{\g,\t} (v_1,w)$. If
$ \a_k := \min_{j \geq 1, j \neq k}$ $| \om^2 k^2 - \lambda_{k,j} | $ 
is attained at 
$j = j(k ) $, i.e. $  \a_k = | \om^2 k^2 - \lambda_{k,j} |$ then, 
$ j \leq 2 k $ 
and therefore, by (\ref{estilb}) and since $1< \tau < 2$, 
we obtain (\ref{akl}).
\end{pf}

\begin{corollary}\label{cor:d-1}
{\bf (Estimate of  $|D|^{-1 \slash 2}$)}
If $ \delta \in \Delta_n^{\g,\t} (v_1,w) \cap [0, \d_0) $ and $\d_0 $ 
is small enough, then 
$ D: W^{(n)} \to W^{(n)} $ is invertible and $ \forall s' \geq 0 $
\be\label{d-1}
\Big| \Big| |D|^{-1/2} h \Big| \Big|_{\s,s'} \leq 
\frac{C}{\sqrt{\gamma}} ||h||_{\s, s'+ \frac{\tau-1}{2}}
\qquad \forall h \in W^{(n)}.
\ee
\end{corollary}

\begin{pf}
Use (\ref{esdk}) and (\ref{akl}). 
\end{pf}

\noindent
{\large {\sc {\bf Step 2: Inversion of $ {\cal L}_n : W^{(n)} \to
W^{(n)} $}}}.  
\\[3mm]
\indent
In order to show the invertibility of 
${\cal L}_n : W^{(n)} \to W^{(n)}$  
it is a convenient devise to write
$$
{\cal L}_n = D - M_1 - M_2 = |D|^{1/2} \Big( U- R_1 - R_2 \Big) |D|^{1/2}  
$$
where 
$$
U := |D|^{-1/2} D|D|^{-1/2} = |D|^{-1} D \qquad \qquad {\rm and} \qquad \qquad
R_i := |D|^{-1/2} M_i |D|^{-1/2}, \quad i= 1, 2.
$$
We shall prove the invertibility of 
$ U -R_1-R_2 : W^{(n)} \to W^{(n)}$ showing that, for $\e $ small enough,  
$ R_1$ and $ R_2 $ are small perturbations of $ U $.

\begin{lemma}\label{lem:U} {\bf (Estimate of $||U^{-1}||$)}
$ U : W^{(n)} \to W^{(n)}$ is an invertible operator and its inverse
$U^{-1}$  satisfies, 
$ \forall s' \geq 0 $, 
\be\label{Ust} 
\Big| \Big| U^{-1}h \Big| \Big|_{\s,s'}= 
|| h ||_{\s,s'} \Big( 1+ O(\e ||a_0||_{H^1}) \Big)  
 \qquad \forall \ h \in W^{(n)}.
\ee
\end{lemma}

\begin{pf}
$ U_k := |D_k|^{-1} D_k : F_k \to F_k $
being  orthogonal for the $\la \ , \ \ra_\e $ scalar product,
it is invertible and $\forall u \in F_k, ||U_k^{-1} u||_\e= ||u||_e$.
Hence, by (\ref{eq:ue1}), there is $C \geq 0$ such that
$$
\forall u \in F_k, \  \  ||U_k^{-1} u||_{H^1}\leq  ||u||_e
(1+C\e ||a_0||_{H^1}) 
$$    
Therefore, $ U = |D|^{-1} D $, being defined by 
$(Uh)_k= U_k h_k$, $ \forall \ |k | \leq L_n $, $U$ is invertible,
$(U^{-1}h)_k=U_k^{-1} h_k$ and (\ref{Ust}) holds.   
\end{pf} 

For proving the smallness of $R_1$ and $R_2$  
we need the following preliminary Lemma. 

\begin{lemma}
\label{bounda}
There are $\mu>0$, $\d_0 > 0 $ and $ C > 0 $ with the following
  property :
if $||v_1||_{0,s} 
\leq 2R$, $[w]_{\s,s} \leq \mu$ and $\d \in [0,\d_0)$, then 
  $||a||_{\s,s+\frac{2(\t-1)}{\b}} \leq C$.
\end{lemma}
\begin{pf}
By the Definition \ref{wss} of $[w]_{\s,s}$  
there are $h_0, h_1,  
\ldots , h_q \in W^{(i)}$ 
and a sequence $(\s_i)_{0 \leq i \leq q}$ with $ \s_i > \s$, such that 
$w= h_0 + h_1 + \ldots + h_q$ and 
\be\label{sum1}
 \sum_{i=0}^q \frac{||h_i||_{\s_i ,s}}{(\s_i - \s)^{  \frac{2(\t-1)}{\b}
}} \leq 2\mu.
\ee
An elementary calculus shows that\footnote{Using that 
$\max_{k \geq 1} k^{\a} 
\exp\{ - (\s_i - \s) k\}  \leq C ( \a ) \slash (\s_i - \s)^{\a}$}
\be\label{sum2}
||h_i||_{\s,s+\frac{2(\t-1)}{\b}} \leq C(\t) \frac{||h_i||_{\s_i,s}}{
(\s_i -\s)^{\frac{2(\t-1)}{\b}}}. 
\ee
Hence, by (\ref{sum1})-(\ref{sum2})
$$
||w||_{\s,s+ \frac{2(\t-1)}{\b}} \leq 
\sum_{i=0}^q ||h_i||_{\s,s+\frac{2(\t-1)}{\b}}
 \leq \sum_{i=0}^q C(\tau)  
\frac{||h_i||_{\s_i,s}}{(\s_i -\s)^{\frac{2(\t-1)}{\b}}} 
\leq 2C(\tau)\mu
$$ 
and by Lemma \ref{vdue},
provided $\mu$ is small enough, $v_2(\d,v_1,w)$ is well defined and  
$||v_2(\d,v_1,w)||_{\s,s+\frac{2(\t-1)}{\b}} 
\leq C' $. 
Hence, by the algebra property of the norm $|| \ ||_{\s,s+\frac{2(\t-1)}{\b}}$
and the analyticity of $f$, $||a||_{\s,s+\frac{2(\t-1)}{\b}}$ 
is bounded by some 
constant, provided $\delta_0$ has been chosen small enough.  
\end{pf}

The ``smallness'' of $ R_2 : W^{(n)} \to W^{(n)} $ 
is just a consequence of the regularizing 
property of $ \partial_w v_2: X_{\s,s} \to X_{\s,s+2} $ proved in 
Lemma \ref{vdue} and  Lemma \ref{bounda}.

\begin{lemma} \label{lem:R2}
{\bf (Estimate of $R_2$)}
Under the hypotheses of ($P3$), there exists 
a constant $ C >0 $ depending on $ \mu $ such that
$$
\Big| \Big| R_2 h \Big| \Big|_{\s, s + \frac{\t-1}{2}} \leq 
C \frac{| \e |}{\g} \  ||h||_{\s,s + \frac{\t-1}{2}} 
\qquad \forall h \in W^{(n)}.
$$
\end{lemma} 

\begin{pf}
Using (\ref{d-1})
and the regularizing estimates 
$||\partial_w v_2 [u]||_{\s,s+2} \leq C ||u||_{\s,s} $ of Lemma
\ref{vdue} we get
\begin{eqnarray*}
\Big| \Big| R_2 h \Big| \Big|_{\s,s + \frac{\t-1}{2}} 
& \leq & \frac{C}{\sqrt{\gamma}} \Big| \Big| M_2 |D|^{-1/2} h
\Big| \Big|_{\s,s + \t -1 } 
= C \frac{| \e |}{\sqrt{\gamma}} 
\ \Big| \Big| P_n \Pi_W \Big( a \  \partial_w v_2 \Big[ |D|^{-1/2}  
h \Big] \Big) \Big| \Big|_{\s,s+ \t -1} \\
&\leq & C \frac{| \e |}{\sqrt{\gamma}}
 \ || a||_{\s,s+ \t -1} \Big| \Big| \partial_w v_2 \Big[
|D|^{-1/2} h \Big] \Big| \Big|_{\s,s+ \t -1} \\
&\leq &   C' \frac{| \e |}{\sqrt{\gamma}}
 \ || a||_{\s,s+ \t -1} \Big| \Big| \partial_w v_2 \Big[
|D|^{-1/2} h \Big] \Big| \Big|_{\s,s+ 2}\\
&\leq& C \frac{| \e |}{\sqrt{\gamma}} 
\ || a||_{\s,s+ \t -1} \Big| \Big| 
|D|^{-1/2}h \Big| \Big|_{\s,s } 
\leq C' \frac{| \e |}{\gamma}  \ ||h||_{\s,s+ \frac{\t-1}{2}} 
\end{eqnarray*} 
since $ 1 < \t < 2 $ and by Lemma \ref{bounda},
$ || a||_{\s,s + \t -1} \leq || a||_{\s,s+\frac{2(\t-1)}{\b}} 
\leq C $.
\end{pf}

The estimate of the ``off-diagonal'' operator $ R_1 : W^{(n)} \to W^{(n)}$ 
requires, on the contrary, a careful analysis
of the ``{\it small divisors}'' and the
use of the ``first order Melnikov non-resonance condition''
$ \d \in $ $ \Delta_n^{\g,\t} (v_1,w)$,
see Definition \ref{Deltap}.
For clarity, we enounce such property
separately.

\begin{lemma} {\bf (Analysis of the Small Divisors)}\label{sme1}
Let $ \d \in \Delta_n^{\g,\t} (v_1,w) \cap [0, \d_0 )$, with $\d_0 $ small.
There exists $ C > 0 $ such that, $ \forall l \neq k $, 
\be\label{sme}
\frac{1}{\al_k \al_l} \leq 
C \frac{|k-l|^{2 \frac{\tau-1}{\beta}}}{\gamma^2 
{| \e |^{\tau-1}} } \qquad \quad {\rm where} \qquad \quad
\b := \frac{2-\tau}{\tau}.
\ee
\end{lemma}

\begin{pf}
To obtain (\ref{sme})
we distinguish different cases. 
\begin{itemize}
\item {\sc First case:} $ | k - l | \geq [\max(|k|,|l|)]^{\b}$. Then 
 $ (\a_k \a_l )^{-1} \leq C |k-l|^{2 \frac{\tau-1}{\b}} \slash \g^2 $. 
\end{itemize}
Indeed we estimate both $ \a_k $, $ \a_l $ 
with the lower bound (\ref{akl}), $ \a_k \geq C \gamma \slash 
|k|^{\tau-1}$, 
$ \a_l \geq C \gamma \slash |l|^{\tau-1} $ and therefore
$$
\frac{1}{\a_k \a_l} \leq C \frac{|k|^{\tau-1} |l|^{\tau-1}}{\gamma^2} 
\leq C \frac{[\max (|k|, |l| )]^{2(\t-1)}}{\gamma^2} \leq C \frac{|k-l|^{2
\frac{\tau-1}{\b}}}{\g^2}.
$$
\begin{itemize}
\item {\sc Second case:}  $| k - l | < [\max(|k|,|l|)]^{\b} $ and 
($ | k | \leq 1/3 |\e |$ or 
$ | l | \leq 1/3 |\e |) $. Then
$ ( \a_k \a_l)^{-1} \leq C / \g $. 
\end{itemize}
Note that, in this case, $ \sign (l) = \sign (k) $ and,  
to fix the ideas, we assume  in the sequel that $ l, k \geq 0 $ (the 
estimate for $k,l < 0 $ is the same, since $ \a_k \a_l = \a_{-k} \a_{-l} $).

Suppose, for example, that $ 0 \leq  k \leq 1/3 |\e |$. 
We claim that if 
$\e $ is small enough,
then $ \a_k \geq (k+1) \slash 8 $. 
Indeed, $ \forall j \neq k $, 
$$
| \om k - j | = 
\Big| \om k - k + k -j \Big| \geq |k -j | - | \om - 1 | \ | k | 
\geq 1 - 2 |\e | \ k \geq \frac{1}{3}.
$$
Therefore $ \forall 1 \leq k < 1 \slash 3 | \e | $, 
$ \forall j \neq k $, $ j \geq 1 $, 
$ |\om^2 k^2 - j^2 | = $ $|\om k - j | \ |\om k + j | \geq$ 
$(\om k +1) \slash 3 \geq (k+1) \slash 6 $ and so  
$$ 
\a_k := \min_{j \geq 1, k \neq j } \Big| \om^2 k^2 - \lambda_{k,j} \Big| 
=  \min_{j \geq 1, k \neq j }\Big| \om^2 k^2 - j^2 - \e M(\d, v_1 ,w ) + O \Big( 
\frac{\e ||a_0 ||_{H^1} }{j} \Big) \Big| \geq 
\frac{k+1}{6} - |\e | \ C \geq \frac{k+1}{8}.
$$
Next, we estimate $ \a_l $. If also $0 \leq l \leq 1\slash 3 | \e | $
then $\a_l \geq 1 \slash 8 $ and therefore 
$(\a_k \a_l)^{-1} \leq 64$. Otherwise, if $ l > 1\slash 3 | \e | $,
we estimate $ \a_l $  with the lower bound (\ref{akl}) and so, 
since\footnote{Indeed $ | k - l | \leq (\max \{k,l \})^{\b} $
and so or $ l \leq k $ or $ l \leq k + l^{\b} $ and so, $l\slash 2 < k$,
since $l \geq 1 \slash 3 |\e |$.}  
$ l \leq 2k $ and $1 < \t < 2 $   
$$
\frac{1}{\a_k \a_l } \leq C \frac{l^{\tau-1}}{ k \gamma } \leq 
\frac{C'}{ k^{2-\tau} \gamma} \leq \frac{C'}{ \gamma}.
$$
In the remaining cases we consider $|k-l| < [\max(|k|,|l|)]^{\b}$
and both $ |k| , |l| > 1/3 | \e | $. 
We have to distinguish two sub-cases. 
For this, 
$ \forall k \in {\bf Z}$, 
let $ j = j(k) \geq 1 $ be the unique integer such that 
$ \a_k := \min_{n \neq |k|} | \om^2 k^2 - \lambda_{k,n}| = $
$ |\om^2 k^2 - \lambda_{k,j} |$. 
Analogously let $i = i(k) \geq 1 $ be the unique integer such that 
$ \a_l = $ 
$ | \om^2 l^2 - \lambda_{l,i} | $.

\begin{itemize}
\item {\sc Third case:} $0 < |k-l| < [\max(|k|,|l|)]^{\b} $, 
$ |k|,|l|  > 1/3 |\e |$ and  
$|k -l| =  |j - i|$. Then  $(\a_k \a_l)^{-1} 
\leq C/ \g | \e |^{\tau-1}.$
\end{itemize}

Indeed  
$ | (\om k -j) - (\om l -i) |= $ 
$| \om (k-l) - ( j - i )| =$ $|\om -1 | | k-l| \geq$ $ | \e | \slash 2 $
and therefore or $|\om k -j| \geq |\e |  /4$ or $|\om l -i| \geq |\e |/4$. 
It follows that 
$ |\om^2 k^2 - j^2 | = $ $ | \om k - j | \ | \om k + j | \geq $ 
$| \e | \om k \slash 2 \geq $$ |\e | k \slash 3 $
and so, for $ \e $ small enough, 
$ | \alpha_k | \geq |\e |k / 4$, or $| \alpha_l | \geq |\e |l /4 $.
Hence, since $ l \leq 2 k $ and $k > 1 \slash 3 $,    
$$
\frac{1}{\a_k \a_l } \leq C \frac{ l^{\tau-1} }{ \gamma |\e | k} \leq
\frac{C}{ \gamma k^{2-\tau} | \e |} \leq \frac{C}{ \gamma | \e |^{\tau-1}}
$$
\begin{itemize}
\item {\sc Fourth case:} 
$0 < |k-l| < [\max(k,l)]^{\b} $, 
$ k , l > 1/3 |\e |$ and 
$ |k -l|  \neq  |j - i|$. Then  $(\a_k \a_l)^{-1} \leq C/\g^2$.
\end{itemize}
Using that $\om $ is $\g-\tau$-Diophantine, that 
$ |k-l| < [\max(k,l)]^{\b}  $
and so $l \geq 2 k $,  
$$
\Big| (\om k - j) - ( \om l - i ) \Big| 
= \Big| \om ( k - l ) - ( j - i ) \Big| \geq \frac{\g}{|k-l|^{\tau} }
\geq \frac{\gamma}{ [\max (k,l)]^{\b \tau}} \geq
c \frac{\gamma}{ k^{\b \tau} } 
$$
so that or $ |\om k - j | \geq c \gamma / 2 k^{\b \tau} $ or  
$ |\om l - i | \geq c \gamma / 2 l^{\b \tau} $. 
Therefore
$ | \om^2 k^2 - j^2 | \geq  C \gamma k^{1 -\b \tau} $ 
and, for $\e $ small enough,   
$ \a_k \geq C \gamma k^{1 - \b \tau} \slash 2 $.
We estimate $\a_l $ with the worst possible lower bound and so, 
using also  $l \leq 2 k $,
$$
\frac{1}{\a_k \a_l } \leq \frac{C l^{\tau-1}}{\gamma^2 k^{1 - \b \tau}} 
\leq C \frac{k^{\tau -2 +\b \tau}}{\gamma^2 } \leq 
\frac{C}{\gamma^2},
$$
since $ \beta := (2-\tau ) \slash \tau $.
Collecting the estimates of all the previous cases, 
(\ref{sme}) follows.
\end{pf}

\begin{lemma}\label{lem:sd} {\bf (Bound of an off-diagonal operator)}
Assume that $\delta \in \Delta_n^{\g,\t} (v_1,w) \cap [0, \d_0) $ and
let, for some $ s' \geq s $, $ b (t,x) \in X_{\s,s'+ 
\frac{\t-1}{\b}}$ satisfy $ b_0 ( x )=0$, {\it i.e.} 
$\int_0^{2\pi} b(t,x) \ dt \equiv 0,$ $ \forall x \in (0,\pi)$. Defining 
the operator $ T_n : W^{(n)} \to W^{(n)}$ by 
$$
T_n h := |D|^{-1/2} P_n \Pi_W \Big( b(t,x) \ |D|^{-1/2}h \Big),
$$
there is a constant $\wtilde{C}$, independent of $b(x,t)$  and of $p$, 
such that 
$$
\Big| \Big| T_n h \Big| \Big|_{\s,s'} \leq 
\frac{\wtilde{C}}{{|\e |}^{\frac{\t-1}{2}}\gamma} 
||b||_{\s,s'+ \frac{\t-1}{\b}} ||h||_{\s,s'} \qquad \forall h \in W^{(n)}. 
$$ 
\end{lemma}
\begin{pf}
For $h \in W^{(n)}$, we have $(T_n h) (t,x)=\sum_{|k|\leq L_n} (T_n
h)_k (x) \exp
(ikt)$, with
\begin{eqnarray} 
(T_n h)_k & = & |D_k|^{-1/2} \pi_k \Big( b \ 
|D|^{-1/2} h \Big)_k \nonumber \\
&=&  |D_k|^{-1/2} \pi_k  \Big[ 
 \sum_{|l| \leq L_n } b_{k-l} | D_l |^{-1/2} h_l   \Big]
\label{estp1}
\end{eqnarray} 
(since $ {\ov b}_{k-l}(x) = b_{l-k} (x) $,
the linear operator $ T_n $ is represented by a self-adjoint Toepliz matrix in
Mat$ (2L_n \times 2L_n,$ $H^1 ((0,\pi), {\bf C}) ) $ 
 which is zero on the 
diagonal, as $b_0 (x) = 0 $). Abbreviating
$ B_m := ||b_m(x)||_{H^1} $,
we get from (\ref{estp1}) and (\ref{esdk}), using that
$B_0 := ||b_0(x)||_{H^1} = 0$,  
\be\label{tpk}
\Big|\Big|(T_n h)_k \Big|\Big|_{H^1} 
\leq C \sum_{|l| \leq L_n , l \neq k}
\frac{B_{k-l}}{\sqrt{\al_k}
\sqrt{\al_l}} ||h_l||_{H^1}.  
\ee
Hence, by (\ref{sme}), 
\be\label{tpsm}
\Big| \Big| (T_n h)_k \Big| \Big|_{H^1} 
\leq 
\frac{C}{\gamma {| \e |}^{\frac{\tau-1}{2}}} s_k \qquad
{\rm where} \qquad 
s_k:= 
\sum_{|l| \leq L_n}
B_{k-l}|k-l|^{\frac{\t-1}{\b}} ||h_l||_{H^1}.  
\ee
By (\ref{tpsm}),
setting $\wtilde{s}(t)
:=\sum_{|k| \leq L_n} s_k \exp (\ii k t)$,
\begin{eqnarray}\label{eqf}
\Big| \Big| T_n h \Big| \Big|_{\s,s'}^2 & = & 
\sum_{|k| \leq L_n} \exp{(2 \s |k|)} ( k^{2s'}+1 ) 
\Big| \Big| (T_n h)_k \Big| \Big|_{H^1}^2  \nonumber \\
&\leq & \frac{C^2}{\gamma^2 |\e |^{\t-1}}\sum_{|k| \leq L_n} \exp{(2\s |k|)} 
( k^{2s'}+1 ) s_k^2 = 
\frac{C^2}{\gamma^2 |\e |^{\t-1}} || \wtilde{s} ||_{\s,s'}^2 
\end{eqnarray}
It turns out that 
$ \wtilde{s} = P_n (\wtilde{b}\wtilde{c})$ where
$\wtilde{b}(t):=\sum_{l \in {\bf Z}} |l|^{\frac{\t-1}{\b}} B_l \exp
(\ii lt)$ and
$\wtilde{c}(t):=\sum_{|l| \leq L_n} ||h_l||_{H^1} \exp (\ii lt)$.
Therefore, by (\ref{eqf}),
$$
\Big| \Big| T_n h \Big| \Big|_{\s,s'} \leq 
\frac{C}{\gamma {|\e |}^{\frac{\t-1}{2}}}
||\wtilde{b}\wtilde{c}||_{\s,s'} 
\leq  \frac{C}{\gamma {|\e |}^{\frac{\t-1}{2}}}
||\wtilde{b}||_{\s,s'} ||\wtilde{c}||_{\s,s'}
\leq \frac{C}{\gamma  {|\e |}^{\frac{\t-1}{2}}} ||b||_{\s,s'+ 
\frac{\t-1}{\b}} ||h||_{\s,s'}
$$
since $ ||\wtilde{b}||_{\s,s'}\leq ||b||_{\s,s'+ 
\frac{\t-1}{\b}}$ and $||\wtilde{c}||_{\s,s'} = ||h||_{\s,s'} $.
\end{pf}

\begin{lemma}{\bf (Estimate of $R_1$)} \label{R1}
Under the hypotheses of (P3), there exists 
a constant $ C >0 $ depending on $ \mu $ such that
$$
\Big| \Big| R_1 h \Big| \Big|_{\s, s + \frac{\t-1}{2}} 
\leq |\e |^{\frac{3-\t}{2}} \frac{C}{\g} ||h||_{\s,s+\frac{\t-1}{2}} \qquad 
\forall h \in W^{(n)}.
$$
\end{lemma}

\begin{pf} 
Recalling the definition of $ R_1 := |D|^{-1/2}  M_1 |D|^{-1/2} $ and $ M_1$, 
and using Lemma \ref{lem:sd} since $\ov{a}(t,x) $ has zero time-average, 
\begin{eqnarray*}
\Big| \Big| R_1 h \Big| \Big|_{\s,s + \frac{\t-1}{2}} 
& = & \Big| \Big||D|^{-1/2}
 M_1 |D|^{-1/2} h \Big| \Big|_{\s,s + \frac{\t-1}{2}}   
= | \e | \ \Big| \Big| |D|^{-1/2} P_n \Pi_W \Big( \ov{a} \ |D|^{-1/2} h \Big) 
\Big| \Big|_{\s,s + \frac{\t-1}{2}} \\
& \leq & | \e |  \  \frac{ \wtilde{C} }{{|\e |}^{\frac{\t-1}{2}}\gamma} 
|| \ov{a} ||_{\s, s+ \frac{\t-1}{2} + \frac{\t-1}{\b}} || h ||_{\s,s +
\frac{\t-1}{2}} \leq  
| \e |^{\frac{3-\t}{2}}  \  \frac{\wtilde{C}}{\gamma} || \ov{a} ||_{\s, s+
\frac{2(\t -1)}{\b}} 
|| h ||_{\s,s + \frac{\t-1}{2}} \\
& \leq & 
|\e |^{\frac{3-\t}{2}} \frac{C}{\gamma} ||h||_{\s,s + \frac{\t-1}{2}}
\end{eqnarray*} 
since $0 < \b < 1 $ and, by Lemma \ref{bounda},
$ || \ov{a} ||_{\s,s+\frac{2(\t -1)}{\b}} \leq || a||_{\s,s
+\frac{2(\t -1)}{\b}} \leq C $.
\end{pf}

{\sc Proof of property ($P3$) completed.} 
Under the hypothesis of $(P3)$, 
the linear operator $U : W^{(n)} \to W^{(n)} $ 
is invertible by Lemma \ref{lem:U} and, by Lemmas 
\ref{lem:R2} and \ref{R1}, provided $ \d $
is small enough, 
$ || U^{-1} R_1 ||_{\s,s + \frac{\t-1}{2}} $ and 
$|| U^{-1} R_2 ||_{\s,s + \frac{\t-1}{2}} < 1 \slash 4 $. Therefore also
the linear operator $ U - R_1 - R_2: W^{(n)} \to W^{(n)}$ is 
invertible and its inverse satisfies 
\begin{eqnarray} \label{invu}
|| ( U-R_1-R_2)^{-1}h ||_{\s,s + \frac{\t-1}{2}} 
&=& || ( I-U^{-1}R_1-U^{-1}R_2 )^{-1} U^{-1}h ||_{\s,s +
  \frac{\t-1}{2}} \\
&\leq& 2 ||U^{-1}h ||_{\s,s + \frac{\t-1}{2}}
\leq C ||h||_{\s,s +\frac{\t-1}{2}} \qquad \forall h \in W^{(n)}.
\end{eqnarray}
Hence ${\cal L}_n $ is invertible,  ${\cal L}^{-1}_n = |D|^{-1/2} ( U-R_1-R_2 )^{-1} 
|D|^{-1/2}: W^{(n)} \to W^{(n)} $, and by (\ref{d-1}),
(\ref{invu}),
$$
\begin{array}{lll}
|| {\cal L}_n^{-1} h ||_{\s,s} &=& 
|| |D|^{-1/2} ( U-R_1-R_2 )^{-1} 
|D|^{-1/2} h ||_{\s,s} 
\leq \dps \frac{C}{\sqrt{\gamma}} || ( U-R_1-R_2 )^{-1} 
|D|^{-1/2} h ||_{\s,s+ \frac{\t-1}{2}} \\
& \leq & 
\dps \frac{C'}{\sqrt{\gamma}} || |D|^{-1/2} h ||_{\s,s+
\frac{\t-1}{2}} \leq \dps \frac{C''}{\gamma} 
||h||_{\s,s+ \t-1} \leq  \dps
 \frac{C''}{\gamma}  (L_n)^{\t-1} ||h||_{\s,s}
\end{array}
$$
which completes the proof of property $(P3)$.

\section{Solution of the ($Q_1$)-equation}\label{sec:Q1}

Finally, we have to solve the finite dimensional $(Q_1)$-equation 
\be\label{Q1r}
- \Delta v_1 = \Pi_{V_1} {\cal G}(\d, v_1) 
\ee
where 
$$
{\cal G} ( \d, v_1 )(t,x) := 
g \Big( 
\d, x, v_1(t,x) + \wtilde{w} ( \d, v_1 )(t,x)+
v_2 ( \d, v_1, \wtilde{w} ( \d, v_1 ))(t,x) \Big)
$$
and we have to ensure that there are solutions $( \d, v_1) \in B_\infty $  
for a  set\footnote{
It is rather easy to see that for  $\om$ close to $1$
and $\d=[(\om^2-1)/2]^{1/(p-1)}$, 
the set $E_\om=\{ v_1 \in V_1 \ | \ (\d,v_1) \in B_\infty \}$ 
 contains the whole ball $B(2R, V_1)$ only if $\om$ is strongly nonresonant, 
i.e. belongs to the zero-measure  set $ {\cal W}_\g := $ $ \{ \om \in {\bf R} \ | \
|\om l - j | \geq {\gamma \slash  l}, \ 
\forall j \neq l, l \geq 0, j \geq 1 \} $ 
for some $\gamma >0$.  If $\om$ is strongly nonresonant, the existence
of $2\pi/\om$-periodic solutions can be proved (\cite{BB}-\cite{BB1})
for {\sc any} nonlinearity $ f $. For more general frequencies, the
set $E_\om$ has gaps, which makes the analysis of the $(Q1)$ equation
more delicate.} 
of $ \delta $'s of positive. 
\\[1mm]
\subsection{The set ${\cal A}_p$} \label{subsec:Ap}
\indent
By Lemma \ref{q20} the $0^{th}$-order ($Q_1$)-equation 
\be\label{0order}
- \Delta v_1 = \Pi_{V_1} {\cal G}(0, v_1) 
\ee
is the Euler-Lagrange equation of the ``reduced'' functional
$ \Psi_0: B(2R,V_1) \to {\bf R} $ defined in (\ref{psi0})
(note that ${\cal G}(0, v_1) = $
$  a_p (x) (v_1 + v_2 (0,v_1, 0))^p $ and 
see formula (\ref{0Q1})).
\\[1mm]
\indent
For any $ a_p (x) \in H^1((0,\pi), {\bf R})$ 
such that condition (\ref{ap}) is verified, 
$ \Psi_0 $, or the functional  $ \wtilde{\Psi}_0 $
obtained replacing $a_p (x)$ by $ - a_p (x)$,  
possesses, by the Mountain-pass Theorem \cite{AR},   
a non-trivial critical point
$ \ov{v}_1 \in V_1 $ with 
$|| \ov{v}_1 ||_{0,s} \leq R$ ($R$ depending on $a_p$).
More precisely, due to time translation invariance, 
we have a circle of critical points. In fact the
functional $\Psi_0$ is invariant under the action of $\R / 2\pi \Z$ on
$ V_1 $ defined by
$$
(\theta * v_1) (t,x) := v_1 (t-\theta , x).
$$
We shall say that a circle of critical points $[v_1] := \{ \theta * v_1
\ ; \theta \in \R / 2\pi \Z \} $ is non degenerate when 
$ \ker \Psi''_0 (v_1) $ is spanned by $(\partial / \partial \theta)
(\theta * v_1)_{|\theta =0}=\partial_t v_1$. 
\\[1mm]
\indent
We recall that, by Lemma \ref{vq}, condition (\ref{ap}) 
holds for any $ a_p (x) \in {\cal O} \subset 
H^1 (( 0, \pi ), {\bf R})$ where 
$$ 
{\cal O} := \Big\{ a(x) \in H^1 (( 0, \pi ), {\bf R}) \ {\rm such \ that}  \ 
\cases{ a(\pi - x) \neq - a(x), \ {\rm for \ some  \ } x \in [0,\pi],\ 
{\rm if} \ p \  {\rm is \ odd} \cr
a(\pi - x) \neq a(x), \ {\rm for \ some \ } x \in [0,\pi], \ {\rm if} 
\ p  \ {\rm is \ even}} \Big\}. 
$$
We can define similarly the nondegenerate critical circles of $\Phi_0$
in $V$. 
\begin{remark}\label{nondeg}
$[\ov{v}_1] $ is a non-degenerate 
critical circle of $\Psi_0 : V_1 \to {\bf R}$ 
iff $[\ov{v}] = [\ov{v}_1 + v_2 (0, \ov{v}_1, 
0 )] $ is a  non-degenerate critical circle 
of $\Phi_0 :  V \to {\bf R}$, 
i.e. 
iff $\ov{v} \in V $ is, up to time translations, a non-degenerate 
solution  of 
the $0^{th}$-order bifurcation equation (\ref{eq:unpe}). 
Moreover, if $[\ov{v}]$ is a nondegenerate critical circle of $\Phi_0$ 
of critical value $\leq c$ then $[\Pi_{V_1}\ov{v}]$ is a nondegenerate 
critical circle of $\Psi_0$.

In \cite{BP3} it is proved that,
for the nonlinearity $ f ( x, u )  = u^3 $,
the critical points of $ \Phi_0 $ are non-degenerate
up to time translations. 
\end{remark}

\noindent
Let us define
$$
{\cal A}^1_p := \Big\{ a_p (x) \in {\cal O} \subset H^1(0,\pi) 
\ | \ 
\hbox{there is a non-degenerate critical circle }
[\ov{v}] \subset V\backslash \{ 0\}  \ {\rm of} \ 
\Phi_0 \ {\rm or} \ \wtilde{\Phi}_0 
\Big\}, 
$$
where $\wtilde{\Phi}_0$ is obtained from $\Phi_0$ replacing $a_p$ with 
$-a_p$.  

For any $ a_p \in {\cal A}^1_p $ 
the $0^{th}$-order bifurcation equation (\ref{eq:unpe}) 
(or the one
obtained substituting $ - a_p $ for $ a_p $)   
possesses a non-trivial, non-degenerate solution (up to time translations)
$ \ov{v} \in V $.  We can always assume that we have chosen $c$ 
large enough in (\ref{defK0}), so that, by Lemma \ref{q20}, $\ov{v}=
 \ov{v}_1+v_2(0,\ov{v}_1,0)$ for some non degenerate (up to time
 translations) solution $\ov{v}_1 \in B(2R,V_1)$ of equation (\ref{0order}).
By the Implicit function Theorem, there exists 
a $C^\infty$-curve of solutions of the $(Q_1)$-equation (\ref{Q1r}) 
$$
v_1 ( \cdot ) : \ [0, \d_0) \to V_1 \qquad 
{\rm with} \qquad  v_1 (0) = \ov{v}_1.
$$
Let us  define the Cantor-like set  
$$
{\cal C}_{a_p, \ov{v}_1} := \Big\{ \d \in [0,\d_0) \ | \ (\d, v_1 (\d)) \in B_\infty 
\Big\}.  
$$ 
The smoothness of  $ v_1 ( \cdot )$ implies that the Cantor set
${\cal C}_{a_p, \ov{v}_1}$ has full density at the origin, i.e.
satisfies the measure 
estimate  (\ref{meas}) of Theorem \ref{thm:main}-$(i)$.

\begin{proposition}\label{measure} {\bf (Measure estimate of 
$ { \cal C}_{a_p,\ov{v}_1} $)}
$\forall a_p ( x )  \in {\cal A}^1_p $,  
$\lim_{\eta \to 0^{+}}$ ${\rm meas}$ $({\cal C}_{a_p , \ov{v}_1} 
\cap (0, \eta)) \slash \eta = 1$.
\end{proposition}

\begin{pf}
Recall that 
\begin{eqnarray*}
B_\infty :=  \cap_{n \geq 1} B_n & = &
\Big\{ (\d, v_1) \in A_0 \ \ :  \ \
\Big| \om (\d ) l - j - \d^{p-1} \frac{M(\d, v_1, \wtilde{w}(\d, v_1))}{2j} 
\Big| \geq \frac{2 \g }{(l+j)^{\tau}},  \\
& & \ \ 
\Big| \om (\d ) l - j  \Big| \geq \frac{2\g}{(l+j)^{\tau} }, \
\forall l, j \geq \frac{1}{3 \d^{p-1}}, \ l \neq j \ \Big\} 
\end{eqnarray*}
where $\om (\d) = \sqrt{1 + 2 \d^{p-1}}$ (or $\om (\d) = \sqrt{1 - 2 \d^{p-1}}$), 
$ M (\d, v_1, w) $ 
is defined in Definition \ref{def:MV} and  $\wtilde w (\d, v_1) $ 
in Lemma \ref{smoex}.

Let $ 0 < \eta < \d_0 $ . The complementary set of $ {\cal C}_{a_p , \ov{v}_1} $ 
in $(0 , \eta )$ is 
\begin{eqnarray*}
{\cal C}_{a_p , \ov{v}_1}^c 
&:= & \Big\{ \d \in (0, \eta ) \ \Big| \ \
\Big| \om( \d ) l - j - \frac{\d^{p-1} m(\d)}{2j} \Big| < 
\frac{2 \gamma}{(l+j)^{\tau}} \quad {\rm or} \quad  
\Big| \om( \d ) l - j \Big| < 
\frac{2 \gamma}{(l+j)^{\tau}} \\
& &  { \rm for \ some \ }
l, j > \frac{1}{3\d^{p-1}}, \ l \neq j \ \Big\}
\end{eqnarray*}
where
$ m(\d ) := M (\d, v_1 (\d), \wtilde{w} (\d, v_1 (\d )) )$
is a function 
in $ C^\infty ([0, \d_0), {\bf R})$ since $ \d \mapsto v_1 (\d ) $ 
is $ C^\infty $ and
$ \wtilde w (\d, v_1) $ is, by Lemma \ref{smoex}, in
$ C^{\infty}(A_0, W \cap X_{\ov{\s}\slash 2, s})) $. 
This implies, in particular,  
\be\label{boun}
| m(\delta)| + | m'(\delta)| \leq C, \qquad \forall \d \in [0,\delta_0/2]
\ee
for some positive constant $ C $.

We claim that, for any interval 
$[\d_1 \slash 2, \d_1] \subset [ 0, \eta ] \subset [0, \d_0 /2] $
the following measure estimate holds:
\be\label{claimc}
{\rm meas}\Big( {\cal C}_{a_p , \ov{v}_1 }^c \cap \Big[ \frac{\d_1}{2} , \d_1 \Big] 
\Big) \leq K_1 (\tau) \gamma \eta^{ (p-1)(\tau-1)} {\rm meas} 
\Big( \Big[ \frac{\d_1}{2} , \d_1 \Big] \Big)
\ee
for some constant $K_1 ( \tau ) > 0 $.

Before proving (\ref{claimc}) we show how to conclude the proof of
the Lemma.
Writing $ (0,\eta] = $ $ \cup_{n \geq 0} $ $ [\eta \slash 2^{n+1},$
$\eta \slash 2^n ]$ and applying the measure estimate (\ref{claimc}) to any
interval $ [\d_1 \slash 2, \d_1] = [\eta \slash 2^{n+1},
\eta \slash 2^n ]$, we get
$$
{\rm meas} ( {\cal C}_{a_p , \ov{v}_1}^c \cap [0,\eta]) \leq 
K_1 (\tau) \gamma \eta^{ (p-1)(\tau-1)} \eta , 
$$
whence $ \lim_{\eta \to 0^{+}}$ ${\rm meas}$ $({\cal C}_{a_p , \ov{v}_1} 
\cap (0, \eta)) \slash \eta = 1$, proving the Lemma. 
\\[1mm]
We now prove (\ref{claimc}). We have
\be\label{resto}
{\cal C}_{a_p , \ov{v}_1 }^c \bigcap \Big[ \frac{\d_1}{2} , \d_1 \Big] \subset
\bigcup_{(l,j) \in {\cal I}_R } {\cal R}_{l,j}( \d_1 )
\ee
where 
$$
{\cal R}_{l,j}( \d_1 )  := \Big\{ \delta 
\in \Big[ \frac{\d_1}{2} , \d_1 \Big]
\ | \
\Big| \om( \d ) l - j - \frac{\d^{p-1} m(\d)}{2j} \Big| < 
\frac{2 \gamma}{(l+j)^{\tau}} \quad {\rm or} \quad  
\Big| \om( \d ) l - j \Big| < 
\frac{2 \gamma}{(l+j)^{\tau}} \Big\}
$$
and 
$$
I_R := \Big\{ (l,j) \ | \ l, j > \frac{1}{3 \d_1^{p-1}}, 
\  l \neq j, \ \frac{j}{l}
\in [1- c_0 \d_1^{p-1} , 1+ c_0 \d_1^{p-1}] \Big\}
$$
(note indeed that 
$ {\cal R}_{j,l} ( \d_1 ) = $ $ \emptyset $ unless 
$ j \slash l \in [1- c_0 \d_1^{p-1} , 1+ c_0 \d_1^{p-1}]$ for some
constant $ c_0 > 0 $ large enough).

Next, let us prove that 
\be\label{rlj}
{\rm meas} ( {\cal R}_{lj}(\d_1)  ) = O 
\Big(\frac{\gamma}{l^{\tau +1} \delta_1^{p-2}} \Big).
\ee
Define 
$ f_{lj}(\d) := $ $ \om( \d ) l - j - (\d^{p-1} m(\d) \slash 2j )$ 
and $ {\cal S}_{j,l}( \delta_1 ) := $ $ \{ \delta \in [ 
\delta_1 \slash 2, \delta_1 ] \ : \ 
|f_{l,j} (\delta )| < 2\gamma \slash (l+j)^\tau \}.$
Provided $ \d_0 $ has been chosen small enough (recall that 
$ j, l \geq  1 \slash 3\d_0^{p-1} $),
$$
f'_{lj} (\delta) = \frac{l (p-1) \d^{p-2}}{\sqrt{1+2 \d^{p-1} }}
- \frac{ (p-1) \d^{p-2} m( \d )}{2j} - \frac{\d^{p-1} m'( \d)}{2j} 
\geq  \frac{(p-1)\delta^{p-2}}{2} \Big( l- \frac{C}{j} \Big) 
\geq \frac{(p-1)\delta^{p-2}l}{4}
$$
and therefore
$f'_{lj} (\delta) \geq (p-1) \delta_1^{p-2} l / 2^p$ for any $\delta
\in [\d_1/2 , \d_1 ]$. This implies
$$
{\rm meas} ({\cal S}_{lj} (\d_1)) \leq  
\frac{4\gamma}{(l+j)^\tau} \times \Big( 
\min_{\d \in [\d_1\slash 2 , \d_1]} f_{lj}'(\d)\Big)^{-1} 
\leq \frac{4\gamma}{(l+j)^\tau} \times 
\frac{2^p}{(p-1) l \d_1^{p-2}} = O 
\Big(\frac{\gamma}{l^{\tau +1} \delta_1^{p-2}} \Big).
$$
Similarly we can prove 
$$
{\rm meas } \Big( \Big\{ \d \in \Big[ \frac{\d_1}{2}, \d_1 \big] \  :  \ 
| \om( \d ) l - j | < \frac{2 \gamma}{(l+j)^{\tau}} \Big\} \Big)  
= O \Big(\frac{\gamma}{l^{\tau +1} \delta_1^{p-2}} \Big)
$$
and the measure estimate (\ref{rlj}) follows. 

Now, by (\ref{resto}), (\ref{rlj}) and since,
for a given $ l $, the number of $ j $ for which $(l,j) \in I_R$ 
is $ O ( \d_1^{p-1} l ) $, 
$$
{\rm meas} \Big({\cal C}_{a_p , \ov{v}_1}^c \cap \Big[ \frac{\d_1}{2}, \d_1 \Big] \Big)
\leq  
\sum_{(l,j) \in I_R} {\rm meas} ({\cal R}_{j,l}(\delta_1) ) \leq 
C \sum_{ l \geq 1 \slash 3 \d_1^{p-1} } \d_1^{p-1} l \times 
\frac{\gamma}{l^{\tau +1} \delta_1^{p-2}} 
\leq K_2 (\tau) \gamma \d_1^{1+ (p-1)(\tau-1)}
$$
whence we obtain (\ref{claimc}) since $ 0 < \d_1 < \eta $.
\end{pf}

Now for $a_p \in {\cal A}^1_p$,  formula 
$$
\wtilde{u}(\d) :=\d \Big[ v_1(\d)+\wtilde{w}(\d, v_1(\d)) + 
v_2 \Big( \d, v_1(\d), \wtilde{w}(\d, v_1(\d))\Big)\Big]
$$
defines a smooth path $\wtilde{u}:[0,\d_0) \to X_{\ov{s}/2,s}$, and 
$\wtilde{u}(\d)=\d u_0 + O(\d^2)$ with $u_0=\ov{v}_1+v_2(0,\ov{v}_1,0)
\in V$. By  Lemma \ref{Binfty}, $\wtilde{u}$ is a solution of 
(\ref{eq:freq}) for $\d \in {\cal C}_{a_p, \ov{v}_1}$ and
conclusions $(i)$ and $(ii)$ 
of Theorem \ref{thm:main} hold
by Proposition \ref{measure}.
\\[2mm] 
Now we can look for $2\pi/(n\om)$ time-periodic solutions of
(\ref{eq:main})
({\it i.e.} $2\pi/n$ time-periodic solutions of
(\ref{eq:freq})) as well. Let
$$
X_{\s,s,n}:= \Big\{ u \in X_{\s,s} \ | \  u \ {\rm is} \ \frac{2\pi}{n} \
{\rm time-periodic}  \Big\}.
$$ 
Replacing $X_{\s,s}$ with $X_{\s,s,n}$, we can develop similarly the
arguments of sections \ref{sec:Q2} and \ref{sec:P}. Define the linear
map  $ {\cal H}_n : V \to V $ by :  for 
$v(t,x) = \eta (t+x) - \eta (t-x) \in V $,
$$
({\cal H}_n v)(t,x) := \eta (n (t+x)) - \eta (n (t-x)) 
$$
and denote by $V_n := {\cal H}_n V $ (resp. $W_n$) the subspace of $V$
(resp. $W$) formed by the
functions $2\pi \slash n $-periodic in $ t $. 

Using the decomposition $X_{\s,s,n}= V_n \oplus W_n$ and introducing an
appropriate finite dimensional subspace $V_{1,n}$ of $V_n$, we
obtain associated $(Q1),(Q2),(P)$-equations (like in (\ref{eqs})),
which can be solved exactly as in the case $ n = 1 $. 

The $0^{th}$-order bifurcation equation is the same (but in $V_n$) and
the corresponding functional is just the restriction of $ \Phi_0 $ to
$V_n$. As $ \Phi_0 $, $\Phi_{0|V_n}$ (or ${\widetilde \Phi}_{0|V_n}$)
possesses ``mountain pass''
critical circles.  Let, for
$n\geq 2$, 
$$
{\cal A}^n_p := \Big\{ a_p (x) \in {\cal O} \subset H^1(0,\pi) 
\ | \ 
\hbox{there is a non-degenerate critical circle }
[\ov{v}] \subset V  \ {\rm of} \
\Phi_{0|V_n} \ {\rm or} \ {\widetilde \Phi}_{0|V_n} \ 
\Big\}, 
$$   
and ${\cal A}_p = \cup_{n=1}^{\infty} {\cal A}^n_p$. 
By the implicit function theorem,  ${\cal A}_p$ is an
open subset of $H^1(0,\pi)$. 

By the arguments of Proposition \ref{measure}, we obtain that if 
$a_p \in {\cal A}_p$, then  conclusions $(i)$ and $(ii)$ of Theorem 
\ref{thm:main} hold (since $2\pi/(n\om)$ time-periodic solutions are
just peculiar $2\pi/\om $ time-periodic solutions).  

This proves the first part of Theorem  \ref{thm:main}.

\subsection{Case  $ f(x,u) = a_3(x) u^3 + O(u^4) $} \label{subsec:u^3}

We assume here that $f(x,u) = a_3(x) u^3 + O(u^4) $, where 
\be\label{media3}
\frac{1}{\pi} \int_0^\pi a_3 (x) dx := \langle a_3 \rangle \neq 0 \ .
\ee
To fix the ideas, we deal with the case  $\la a_3 \ra > 0$.  

Note that assumption (\ref{media3}) also implies 
condition (\ref{ap}), i.e. that there is $ v \in V $ such that
$ \int_\Om a_3(x) v^4 \neq 0 $.  
Indeed  we know that $ \int_\Om a_3(x) v^4 = 0 $, $ \forall v \in V $
iff $ a_3(\pi -x ) = - a_3(x)$, $ \forall x \in [0,\pi] $. But in this  
case $\langle a_3 \rangle = 0 $.

Therefore the $0^{th}$-order bifurcation equation 
is  (\ref{eq:unpe}) (with $p=3$), {\it i. e.} 
the Euler-Lagrange equation of 
\be\label{Phi0u}
\Phi_0 (v) = \frac{||v||_{H^1}^2}{2} - 
\int_\Om a_3(x) \frac{v^4}{4} 
\ee

The functional $ \Phi_n (v) = \Phi_0 ({\cal H}_n v ) $ 
has the following development:
for $ v (t,x) = \eta (t+x) - \eta (t-x) \in V $ we obtain, using that
$ \int_\Om v^4 = \int_\Om ({\cal H}_n v)^4 $,  
$$
\Phi_n (v)  =  2\pi n^2 \int_{{\bf T}} {\dot \eta}^2 (t) dt 
- \langle a_3 \rangle \int_{\Om} \frac{v^4}{4}  \ -  
\int_\Om \Big( a_3(x) - \langle a_3 \rangle \Big) 
\frac{({{\cal H}_n v})^4}{4}. 
$$
Hence
\begin{eqnarray*}
\Phi_n \Big( \frac{\sqrt{2} n}{\sqrt{  \langle a_3 \rangle}} v  \Big) 
& = & \frac{8 \pi n^4}{\la a_3 \ra}\Big[ 
\frac{1}{2} \int_{{\bf T}} {\dot \eta}^2 (s) ds -\frac{1}{8\pi}
 \int_\Om v^4 \ + \frac{1}{8\pi} \int_\Om 
\Big(\frac{a_3(x)}{\la
  a_3\ra } -1 \Big) ({\cal H}_n  v)^4 \ dt \ dx 
\Big] \\
&=& \frac{8\pi n^4}{\la a_3 \ra} [ \Psi(\eta)+ R_n(v) ]
\end{eqnarray*}
where 
$$
\Psi(\eta) := \frac{1}{2} \int_{{\bf T}} {\dot \eta}^2 (s) ds
-   \frac{1}{4} \int_{{\bf T}}  \eta^4 (s) ds -
\frac{3}{8\pi} \Big( \int_{{\bf T}} {\eta}^2 (s) ds \Big)^2,
$$
$$
R_n (v) := \frac{1}{8\pi} \int_\Om b(x) ({\cal H}_n  v)^4 \ dt \ dx, \quad 
b(x) :=\frac{a_3(x)}{\la a_3 \ra}-1.
$$ 

To get that $a_3 \in {\cal A}_3$, 
it is enough to prove that  $ \Psi $ has a non-degenerate critical circle 
and that $R_n$ is a small perturbation for large $n$, more precisely that  
$D^2 R_n \to 0$, $DR_n \to 0$   uniformly on bounded sets as 
$n \to + \infty$.  Then, by the implicit function theorem, for $n$
large enough,  $\Phi_n$
too (hence $\Phi_{0|V_n}$) has a non-degenerate critical circle, 
which implies that $a_3 \in {\cal A}_3$.   

The critical points of $ \Psi $ in 
$ E := $ $\{ \  \eta \in H^1(\T) \ | \ \int_{\T} \eta =0 \ \}$  
are the $ 2\pi $-periodic 
solutions with zero mean value of 
\be
{\ddot \eta} +   
\eta^3 + 3 \langle \eta^2 \rangle \eta  = C, \quad
C \in \R
\ee 
By [BP] it is known that there exists
a solution to this problem (with $ C = 0 $)  which is a
non-degenerate (up to time translations) critical point of $\Psi$
in $ E $. Finally

\begin{lemma} There holds
\be\label{d2n}
|| D R_n ( v ) ||, \ 
|| D^2 R_n ( v ) || \to 0 \qquad {\rm as} \qquad n \to + \infty 
\ee
uniformly for $v$ in bounded sets. 
\end{lemma}

\begin{pf}
We shall  prove the estimate only for $D^2 R_n $. We have
\begin{eqnarray*}
| D^2 R_n ( v )[h,k]| 
& = & \frac{3}{2\pi} \int_\Om b(x) 
({\cal H}_n v )^2 ({\cal H}_n h) \ ({\cal H}_n k)  \\
& = & \frac{3}{2\pi} \int_0^\pi b(x) g( nx ) \, dx 
\end{eqnarray*}
where $g(y)$ is the $ \pi $-periodic function defined by
$$
g(y) := \int_{{\bf T}} (\eta (t+ y) - \eta (t - y))^2 (\beta(t +y) -
\beta(t-y))( \gamma(t+y) - \gamma(t-y)) dt \ , 
$$
$\beta$ and $\gamma$ being associated with $h$ and $k$ as 
$\eta$ is with $v$. 
Developing in Fourier series $g(y) = \sum_{l \in {\bf Z}} g_l \exp (\ii 2 l y )$
we have $g(nx) = \sum_{l \in {\bf Z}} g_l \exp (\ii 2 l n x )$.
Extending
$ b(x) $ to a $\pi$-periodic function, we also write
$ b(x) = \sum_{l \in {\bf Z}} b_l \exp (\ii 2 l x ) $, with $b_0=\la b
\ra =0$. Therefore
\begin{eqnarray*}
|D^2 R_n ( v_n )[h,k]| &=& \frac{3}{2}
\Big| \sum_{l \neq 0} g_l b_{-ln} \Big| \leq 
\frac{3}{2} \Big( \sum_{ l \neq 0 } g_{l}^2 \Big)^{1 \slash 2}
 \Big( \sum_{ l \neq 0 } b_{ln}^2 \Big)^{1 \slash 2}
\\
& \leq & \frac{3}{2}  ||g||_{L^2(0,\pi)} \Big( \sum_{ l \neq 0 } b_{ln}^2 \Big)^{1 \slash 2}\\
& \leq & C || v_0 ||_\infty^2 ||\beta ||_\infty ||\gamma ||_\infty \ 
\Big( \sum_{ l \neq 0 } b_{ln}^2 \Big)^{1 \slash 2}. 
\end{eqnarray*}
Since 
$ ( \sum_{l \neq 0 } b_{ln}^2 )^{1 \slash 2} \to 0 $ as 
$n \to \infty $ it proves (\ref{d2n}). 
With a similar calculus we can prove that $D R_n ( v )\to 0 $
as $ n \to + \infty $.
\end{pf}

This completes the proof of part 1) of Theorem \ref{thm:main}.

\subsection{Case $f(x,u) = a_2 u^2 + O(u^4) $} \label{subsec:u^2}

We now prove part 2) of Theorem \ref{thm:main}. 

In this case condition (\ref{ap}) is violated. It turns out (see
\cite{BB1}) that we must take $\om <1$ ({\it i.e.} $\e<0$), and 
$\d=|\e|^{1/2}$ in the rescaling.  By
the computations of \cite{BB1},   
the $ 0^{th} $-order bifurcation equation is the Euler Lagrange 
equation of the functional $ \Phi_0 : V \to {\bf R} $ defined by  
$$
\Phi_0 (v) = \frac{||v||_{H^1}^2}{2} + \frac{a_2^2}{2} 
\int_\Om v^2 L^{-1} v^2, 
$$
where $ L^{-1} : W \to W $ is the inverse operator of 
$ -\partial_{tt} + \partial_{xx} $.   
\\[1mm]
\indent
We can still use the same arguments to solve the $(Q2)$ and $(P)$
equations.  As explained in subsection \ref{subsec:Ap}, if 
the functional $  \Phi_n : V  \to {\bf R} $, defined by $ \Phi_n (v) = 
\Phi_0 ({\cal H}_n v )$ possesses a nondegenerate critical circle, then
$(i)$ and $(ii)$ hold in Theorem \ref{thm:main} (with
$\om=\sqrt{1-2\d^2}$).  

$\Phi_n $ admits the following development (Lemmae 3.7 and 3.8 in \cite{BB1}): 
for $v(t,x) = \eta (t+x) - \eta (t-x)$ 
\begin{equation}
\Phi_n (v)  =   2\pi n^2 \int_{{\bf T}} {\dot \eta}^2 (t) dt 
- \frac{\pi^2 a_2^2}{12} \Big( \int_{{\bf T}} \eta^2 (t) \ dt \Big)^2 +
\frac{a_2^2}{2n^2}  \Big( \int_\Om v^2 L^{-1} v^2 + \frac{\pi^2}{6}
\Big( \int_{{\bf T}} \eta^2 (t) \ dt \Big)^2 \Big). 
\end{equation}
Hence we can write
\be
\Phi_n \Big( \frac{\sqrt{12}n}{\sqrt{\pi} a_2} v \Big) = \frac{48n^4}{a_2^2} \Big[
\frac{1}{2} \int_{\T} \dot{\eta}^2 (s) \ ds - \frac{1}{4} \Big(   
\int_{\T} {\eta}^2 (s) \ ds  \Big)^2
+ \frac{1}{n^2} R(\eta) \Big] =\frac{48n^4}{a_2^2} \Big[ \Psi (\eta ) + 
\frac{1}{n^2} R(\eta) \Big]
\ee 
where 
\be
\Psi (\eta) =
\frac{1}{2} \int_{\T} \dot{\eta}^2 (s) \ ds - \frac{1}{4} \Big(   
\int_{\T} {\eta}^2 (s) \ ds  \Big)^2
\ee
and
$ R : E \to \R$ is a smooth map defined on  
$E :=\{  \eta \in H^1(\T) \ | \ \int_{\T} \eta =0 \}$.
In order to  prove that $\Phi_n $ has a non-degenerate critical circle for
$n$ large enough, it is enough to prove that:

\begin{lemma} 
$\Psi: E \to {\bf R}$
possesses a  critical point which is non degenerate, up to time translations. 
\end{lemma}

\begin{pf}
The critical points of $ \Psi$ in $E$ are the
$2\pi$-periodic  solutions of zero mean value of the equation 
\be\label{riscanu2}
{\ddot \eta} + 
\Big( \int_{{\bf T}} \eta^2 (t) \ dt \Big) \eta = C, \quad C \in \R. 
\ee
Since $\ddot{\eta}$ and $\eta$ have both zero mean value, any solution
of  (\ref{riscanu2})  must satisfy $C=0$. 
Equation (\ref{riscanu2}) (with $C =0$) 
has a $ 2 \pi$-periodic solution of the form
$ \bar \eta (t)= (1/\sqrt{\pi})\sin t  $. 

We claim that $\bar \eta $ is non-degenerate, up to
time-translations, i.e.
${\rm Ker} \ (D^2 \Psi )( \bar \eta ) = \langle \dot {\bar \eta} \rangle $.
 
The linearized equation
of (\ref{riscanu2}) at $\bar \eta $ is 
$$
 {\ddot h} +  
\Big( \int_{{\bf T}} \bar \eta^2 (t) \ dt \Big) h + 
2 \Big( \int_{{\bf T}} \bar \eta (t) h  (t) \ dt \Big) \bar \eta =c,
$$
and again $c$ must be equal to $0$.
We get
\be\label{equl}
{\ddot h} + h + 
\frac{2}{\pi}
\Big( \int_{{\bf T}} \sin t \ h  (t) \ dt \Big) \sin t = 0 .
\ee
Developing in time-Fourier series 
$$
h(t) = \sum_{k \geq 1} \Big( a_k \sin kt 
+ b_k  \cos kt \Big)
$$
we find out that any solution of the linearized equation (\ref{equl})
satisfies
$$
-k^2 b_k + b_k = 0, \ \ \forall k \geq 1,  \qquad
- k^2 a_k + a_k = 0, \ \ \forall k \geq 2, \qquad 
a_1 = 0 \ 
$$
and 
therefore $ h \in \langle \cos t \rangle 
= \langle \dot {\bar \eta} \rangle  $.
\end{pf}

By the previous Lemma, for $n$ large enough there is  a non-degenerate 
critical circle  
 of $ \Phi_n $ in $ V $. This completes the proof of
Theorem \ref{thm:main}.

\section{Appendix}

\begin{lemma} \label{vq}
If $q$ is an even integer, then  
$$
\int_\Omega a (x) v^q(t,x) \ dt \ dx  = 0, \ \forall v \in V \ 
\Longleftrightarrow \ \ \Big\{ a(\pi - x) = - a(x), \ \forall 
x \in [0,\pi] \Big\}.
$$ 
If  $q \geq 3$ is an odd integer, then
$$ 
\int_\Omega a (x) v^q (t,x) \ dt \ dx  = 0, \ \forall v \in V \ 
\Longleftrightarrow   \ \   \Big\{ a(\pi - x) = a(x), \ \forall 
x \in [0,\pi] \Big\}.
$$ 
\end{lemma}

\begin{pf} 
We first assume that $q=2s$ is even. If $ a ( \pi - x ) = - a ( x ) $ 
$\forall x \in (0,\pi)$, then, for all $ v \in V$, 
\begin{eqnarray*}
\int_\Om  a(x) v^{2s}(t,x) \ dt \ dx &=& \int_\Om  a(\pi-x)
v^{2s}(t,\pi-x) \ dt \ dx \\
&=& \int_\Om  -a(x) (-v(t+\pi,x))^{2s} \ dt \ dx \\
&=& -\int_\Om  a(x) v^{2s}(t,x) \ dt \ dx
\end{eqnarray*}
and so $\int_\Om  a(x) v^{2s}(t,x) \ dt \ dx=0$. \\[2mm]
Now assume that $\Sigma (v):= \int_\Om  a(x) v^{2s}(t,x) \ dt
\ dx=0$ $\forall v \in V $. Writing that $D^{2s} \Sigma =0$,
we get
$$
\int_\Om a(x) v_1(t,x) \ldots v_{2s}(t,x) \ dt \ dx=0, \qquad
\forall (v_1, \ldots , v_{2s} ) \in V^{2s}.
$$
Choosing  $v_{2s}(t,x)=v_{2s-1}(t,x)= \cos{(lt)} \sin (lx)$, we obtain 
$$
\frac{1}{4} \int_\Om a(x)
v_1(t,x) \ldots v_{2(s-1)} (t,x)(\cos(2lt)+1)(1-\cos(2lx)) \ dt \ dx =0
$$
Taking limits as $l\to \infty$, there results
$  \int_\Om a(x) v_1(t,x) \ldots
v_{2(s-1)}(t,x) \ dt \ dx=0$
$\forall (v_1,\ldots,v_{2(s-1)}) \in V^{2(s-1)}$. 
Iterating this operation, we finally get 
$$
\forall (v_1,v_2) \in V^2 \  \  \int_\Om a(x) v_1(t,x)
v_2(t,x) \ dt \ dx=0, \quad {\rm and} \quad \int_0^\pi a(x) \ dx=0.
$$ 
Choosing $v_1(t,x)=v_2(t,x)=\cos(lt) \sin(lx)$ in the first equality,
we derive that $\int_0^\pi a(x) \sin^2 (lx) \ dx =0$. Hence
$$ 
\forall l\in \N  \ \ \int_0^\pi a(x) \cos(2lx) \ dx=0. 
$$
This implies that $a$  is orthogonal in $L^2(0,\pi)$ to 
$$
F= \Big\{ b \in L^2(0,\pi) \ | \ b(\pi-x)=b(x) \ {\rm a.e.} \Big\}. 
$$
Hence $a(\pi-x)=-a(x)$ a.e., and, since $a$ is continuous, the
identity holds everywhere. 
\\[3mm]
We next assume that $q=2s+1$ is odd, $q\geq 3$. 
The first implication is derived
in a similar way. Now assume that $\int_\Om a(x)
v^q (t,x) \ dt \ dx =0 \  \  \forall v \in V$. We can prove exactly as
in the first part that 
$$
\forall (v_1,v_2,v_3) \in V^3 \ \ 
\int_\Om a(x) v_1(t,x) v_2(t,x)
v_3(t,x) \ dt \ dx=0. 
$$
Choosing $v_1(t,x)=\cos(l_1 t) \sin(l_1 x)$, 
$v_2(t,x)=\cos(l_2 t)
\sin(l_2 x)$, $v_3(t,x)=\cos((l_1+l_2) t) 
\sin((l_1+l_2) x)$ and using
the fact that $\int_0^{2\pi} \cos (l_1t) 
\cos (l_2t) \cos ((l_1+l_2)t)
\ dt \neq 0$, we obtain 
\be \label{a(x)1} 
\begin{array}{l}
\dps \int_0^\pi a(x) \Big[ \sin^2 (l_1 x) 
\sin (l_2 x) \cos (l_2 x) + \sin^2
(l_2 x) \sin (l_1 x) \cos (l_1 x) \Big] \ dx = \\
\dps \int_0^\pi a(x) \sin (l_1 x) 
\sin (l_2 x) \sin \Big( (l_1+l_2)x \Big) \ dx =0.
\end{array}
\ee
Letting $l_2$ go to infinity and taking limits,
(\ref{a(x)1}) yields $\int_0^\pi (1/2) a(x) \sin(l_1x) \cos (l_1 x) \
dx =0$. Hence 
$$
\forall l>0 \  \  \int_0^\pi a(x) \sin (2lx) =0.
$$
This implies that, in $L^2(0,\pi)$, $a$
is orthogonal to $G=\{ b \in L^2(0,\pi) \ | 
\ b(\pi-x)=-b(x) \ {\rm a.e.} \}$. Hence 
$a(\pi-x)=a(x) \ \forall x \in (0,\pi)$. 
\end{pf}
\\[2mm]
\noindent
{\bf Proof of Lemma \ref{interm}. }
Let $K_k(\e)=S_k^{-1}(\e)$ be the selfadjoint compact 
operator of $F_k$ defined by
$$
\la K_k(\e ) u , v \ra_\e = (u,v)_{L^2}, \qquad 
\forall u,v \in F_k 
$$
(in other words $ K_k(\e ) u  $ is the unique weak solution $z\in F_k$
of $ S_k z := u $).

Note that $ K_k(\e)$ is a positive operator, 
i.e. $ \la K_k(\e)u,u \ra_\e >0 $, $ \forall u \neq 0 $, and that $ K_k (\e) $ 
is also selfadjoint for the $L^2$-scalar product.

Let $F_{k,r}=\{u \in F_k \ | \ u((0,\pi)) \subset \R \}$. 
We have $S_k (F_{k,r} \subset F_{k,r}$, and by the spectral theory of
compact selfadjoint operators in Hilbert spaces, 
there is a $\la \ , \  \ra_\e$-orthonormal basis
$(v_{k,j})_{j\geq 1, j\neq k}$ of $F_k$, such that 
$v_{k,j} \in F_{k,r}$ is an eigenvector of $K_k (\e)$ associated to a
positive eigenvalue $\nu_{k,j}(\e)$, the sequence $(\nu_{k,j}(\e))_j$ is
non-increasing and tends to $0$ as $j\to +\infty$.
Each $v_{k,j}(\e)$ belongs to $ D(S_k) $ and is an eigenvector of $S_k$
with associated eigenvalue $\lambda_{k,j}(\e)= 1/\nu_{k,j}(\e)$.
$(\lambda_{k,j}(\e))_{j\geq 1}$ is a sequence non decreasing 
and tending to $ + \infty$ as $j \to + \infty$. 

The map $\e \mapsto K_k(\e)\in {\cal L}(F_k,F_k)$ is differentiable
and $K'_k (\e)=-K_k(\e) M K_k(\e)$, where
$ Mu := \pi_k (a_0 u)$. 

For $u=\sum_{j\neq k} \alpha_j v_{k,j}(\e) \in F_k$, 
$$
\la u,u\ra_\e=\sum_{j\neq k} |\alpha_j|^2 \quad {\rm and} 
\quad (u,u)_{L^2}=\sum_{j\neq k} \frac{|\alpha_j|^2}{\lambda_{k,j}(\e)}.
$$
As a consequence,
\be \label{eigen}
\lambda_{k,j}(\e)=\min \Big\{ \max_{u \in F , ||u||_{L^2}=1 } \la u,u\ra_\e
 \ ; \ F \  \hbox{subspace of $F_k$ of dimension} \  j \ ({\rm if} \ j<
  k) \ , \ j - 1 \ ({\rm if} \ j>k) \Big\}.
\ee
It is clear by inspection that $\lambda_{k,j}(0)=j^2$ and that we can
choose $v_{k,j}(0)=\sqrt{2/\pi} \sin (jx)/j$. Hence, by (\ref{eigen}), 
$|\lambda_{k,j}(\e)-j^2| \leq |\e | \ ||a_0||_\infty <1$, 
from which we derive
\be \label{eigen2}
\forall l \neq j \  \ |\lambda_{k,l}(\e) - \lambda_{k,j}(\e)| \geq
(l+j)-2 \geq 2\min(l,j)-1  \ \ (\geq 1 ). 
\ee
In particular, the eigenvalues $\lambda_{k,j}(\e)$ ($\nu_{k,j}(\e)$)
are simple. By the variational characterization (\ref{eigen}) we also 
see that $\lambda_{k,j}(\e)$ depends continuously on $\e$, and we can assume
without loss of generality that $\e \mapsto v_{k,j}(\e)$ is a
continuous map to $F_k$. 

Let $\varphi_{k,j}(\e) := $ $  \sqrt{\lambda_{k,j}(\e)} v_{k,j}(\e) $.
$(\varphi_{k,j}(\e))_{j\neq k}$ is a $L^2$-orthogonal family
in $F_k$ and
$$
\forall \e \quad \left\{
\begin{array}{l}
K_k(\e) \varphi_{k,j}(\e) = \nu_{k,j}(\e) \varphi_{k,j}(\e) \\
(\varphi_{k,j}(\e) , \varphi_{k,j}(\e))_{L^2} =1
\end{array}
\right.
$$
We observe that the $L^2$-orthogonality w.r.t. $\varphi_{k,j}(\e)$ is
equivalent to the $\la \ , \ \ra_\e$-orthogonality
w.r.t. $\varphi_{k,j}(\e)$, and that $E_{k,j}(\e):=
[\varphi_{k,j}(\e)]^{\bot}$ is invariant under $K_k(\e)$. Using that 
$L_{k,j}:=(K_k (\e)-\nu_{k,j}(\e) I)_{|E_{k,j}(\e)}$ is invertible, 
it is easy to
derive from the Implicit Function Theorem that the maps $(\e \mapsto
\nu_{k,j}(\e))$ and $(\e \mapsto \varphi_{k,j}(\e))$ are
differentiable. 

Denoting by $P$ the orthogonal projector onto $E_{k,j}(\e)$,
we have $$
\varphi'_{k,j}(\e)=L^{-1} (-P K_k'(\e) \varphi_{k,j}(\e))=
L^{-1} (P K_k M K_k  \varphi_{k,j}(\e))
=\nu_{k,j}(\e) L^{-1} K_k P  M  \varphi_{k,j}(\e) , $$
\be \label{eigen3}  \begin{array}{lll}
\nu'_{k,j}(\e)&=& \Big( K'_k (\e)  \varphi_{k,j}(\e), 
\varphi_{k,j}(\e) \Big)_{L^2}
= - \Big( K_kMK_k  \varphi_{k,j}(\e), \varphi_{k,j}(\e) \Big)_{L^2}\\
&=&  - \Big( MK_k  \varphi_{k,j}(\e), K_k\varphi_{k,j}(\e) \Big)_{L^2}
= -\nu^2_{k,j}(\e) \Big( M  \varphi_{k,j}(\e), \varphi_{k,j}(\e) \Big)_{L^2}.  
\end{array}\ee
We have
$$ \nu_{k,j} L^{-1}K_k \Big( \sum_{l \neq j} \alpha_l v_{k,l} \Big)=
\sum_{l\neq j} \frac{\nu_{k,j} \nu_{k,l}}{\nu_{k,l}-\nu_{k,j}}
\alpha_l v_{k,l} = \sum_{l\neq j} \frac{\alpha_l}{\lambda_{k,j}
-\lambda_{k,l}} v_{k,l}. 
$$
Hence, by (\ref{eigen2}), $|\nu_{k,j} L^{-1} K_k P u |_{L^2} \leq
|u|_{L^2}/j$.  We obtain $| \varphi'_{k,j}(\e)|_{L_2}=O(|a_0|_\infty /
j)$. Hence
$$
\Big| \varphi_{k,j}(\e)-\sqrt{\frac{2}{\pi}} \sin (jx)\Big|_{L^2} = 
O \Big( \frac{\e |a_0|_\infty}{j} \Big).
$$
Hence, by (\ref{eigen3}), 
$$
\begin{array}{lll}
\lambda'_{k,j}(\e) &=& \Big( M \varphi_{k,j}(\e), 
\varphi_{k,j}(\e) \Big)_{L^2} =
\dps \int_0^\pi a_0 (x) (\varphi_{k,j})^2 \ dx \\
&=& \dps \frac{2}{\pi} \int_0^\pi a_0(x) (\sin(jx))^2 \ dx + O\Big(
 \frac{\e |a_0|^2_\infty}{j}\Big) 
\end{array}
$$
Writing $\sin^2(jx)=(1-\cos(2jx))/2$, and since $\int_0^\pi
a_0(x) \cos(2jx) \ dx =- \int_0^\pi (a_0)_x (x) \sin(2jx) / 2j \ dx$,
we get
$$
\lambda'_{k,j}(\e)= \frac{1}{\pi} \int_0^\pi a_0(x) \ dx +
O\Big(\frac{ ||a_0||_{H^1}}{j} \Big)= M(\d, v_1, w ) + 
O\Big(\frac{ ||a_0||_{H^1}}{j} \Big).
$$
Hence 
$ \lambda_{k,j}(\e)= $ $j^2 +$ $\e M(\delta,v_1,w) +$
$O (\e ||a_0||_{H^1} \slash j ) $, which is the first estimate
in (\ref{expa}).


\begin{thebibliography}{10}


\bibitem{AR} A. Ambrosetti, P. Rabinowitz, 
{\it Dual Variational Methods in Critical Point Theory and Applications}, 
Journ. Func. Anal, 14, 349-381, 1973.


\bibitem{BP1} D. Bambusi, S. Paleari, 
{\it Families of periodic solutions 
of resonant PDEs}, J. Nonlinear Sci., 11,  69-87, 2001.

\bibitem{BP3} D. Bambusi, S. Cacciatori, S. Paleari,
{\it Normal form and exponential stability 
for some nonlinear string equations}, 
Z. Angew. Math. Phys, 52, 6, 1033-1052, 2001. 

\bibitem{BB} M. Berti, P. Bolle, {\it Periodic solutions of
nonlinear wave equations with general nonlinearities}, 
Comm. Math. Phys. Vol. 243, 2, pp. 315-328, 2003.

\bibitem{BB1} M. Berti, P. Bolle, 
{\it Multiplicity of periodic solutions of nonlinear 
wave equations}, Nonlinear Analysis, 56/7, pp. 1011-1046, 2004.

\bibitem{BB0} M. Berti, P. Bolle, {\it Bifurcation 
of free vibrations for completely resonant wave equations}, 
Boll. Unione Mat. Italiana (8) 7-B, 519-528, 2004.

\bibitem{Bo1} J. Bourgain, {\it  Construction of quasi-periodic solutions 
for Hamiltonian perturbations of linear equations and applications 
to nonlinear PDE}, Internat. Math. Res. Notices  1994,  no. 11.

\bibitem{B3} J. Bourgain, {\it Quasi-periodic solutions of Hamiltonian 
perturbations of $2D$ linear Schr\"odinger equations}, 
Ann. of Math., 148,  363-439, 1998.

\bibitem{B2} J. Bourgain, {\it Periodic solutions of nonlinear wave 
equations},  Harmonic analysis and partial differential equations, 
69--97, Chicago Lectures in Math., Univ. Chicago Press, 1999. 

\bibitem{CY} L. Chierchia, J. You, {\it KAM tori for 1D 
nonlinear wave equations with periodic boundary 
conditions}, Comm. Math. Phys. 211, 2000, no. 2, 497--525. 



\bibitem{CW} W. Craig, E.Wayne, {\it Newton's method and periodic solutions 
of nonlinear wave equation}, Comm. Pure and Appl. Math, 
vol. XLVI, 1409-1498, 1993.

\bibitem{CW1} W. Craig, E.Wayne, {\it 
Nonlinear waves and the $1:1:2$ resonance}, 
Singular limits of dispersive waves, 297--313,
NATO Adv. Sci. Inst. Ser. B Phys., 320, Plenum, New York, 1994.  


\bibitem{FR} E. R. Fadell, P. Rabinowitz, {\it Generalized 
cohomological index theories for the group actions 
with an application to bifurcation questions  for Hamiltonian
systems}, Inv. Math. 45, 139-174, 1978.

\bibitem{FS} J. Fr\"{o}hlich, T. Spencer, {\it Absence of diffusion 
in the Anderson tight binding model for large disorder or low energy},
Comm. Math. Phys.  88 , 1983, 2, 151-184. 

\bibitem{GM} 
G. Gentile, V. Mastropietro, {\it
Construction of periodic solutions of the nonlinear wave equation 
with Dirichlet boundary conditions by the Lindstedt series 
method}, to appear in J.  Math. Pures  Appl.


\bibitem{GMP} 
G. Gentile, V. Mastropietro, M. Procesi {\it
Periodic solutions for completely resonant nonlinear wave equations}, 
to appear on Comm. Math. Phys.

\bibitem{IPT} 
G. Iooss, P. Plotnikov, J.F. Toland. {\it 
The standing wave problem on infinite depth},
C. R. Acad. Sci. Paris, S. I ,338 (5), 425-431, 2004.

\bibitem{IPT1}
G. Iooss, P. Plotnikov, J.F. Toland, {\it 
Standing waves on an infinitely deep perfect fluid under gravity} 
preprint 2004.

\bibitem{LS} B. V. Lidskij, E.I. Shulman, 
{\it Periodic solutions of the equation 
$u_{tt}-u_{xx}+u^3=0$}, 
Funct. Anal. Appl. 22 (1988), 332--333.


\bibitem{KP} S. Kuksin, J. P\"oschel, {\it 
Invariant Cantor manifolds of quasi-periodic oscillations 
for a nonlinear Schr\"{o}dinger equation},  
Ann. of Math, 2,  143,  1996,  no. 1, 149-179.

\bibitem{Mo} J. Moser, {\it Periodic orbits near an Equilibrium 
and a Theorem by Alan Weinstein}, Comm.  Pure  Appl. Math.,
vol. XXIX, 1976.


\bibitem{Po2} J. P\"oschel,
{\it A KAM-Theorem for some nonlinear PDEs},  Ann. Scuola Norm. Sup.
Pisa, Cl. Sci., 23, 1996, 119-148.
\smallskip

\bibitem{Po3} J. P\"oschel, {\it Quasi-periodic solutions for 
a nonlinear wave equation}, Comment. Math. Helv.,  71,  1996, 
no. 2, 269-296.
\smallskip

\bibitem{W1} E. Wayne, {\it Periodic and quasi-periodic solutions 
of nonlinear wave equations via KAM theory},
Commun. Math. Phys. 127, No.3, 479-528, 1990.

\bibitem{We} A. Weinstein, {\it Normal modes for 
Nonlinear Hamiltonian Systems}, Inv. Math, 20, 47-57, 1973. 

\end{thebibliography}
\end{document}